\documentclass[journal]{IEEEtran}
\usepackage{amsmath,amssymb,amsfonts}

\usepackage{epsfig}
\usepackage{subfigure}

\usepackage{threeparttable}
\usepackage{booktabs}
\usepackage{url}
\usepackage{multirow}

\usepackage{fixmath}
\usepackage{graphicx}

\usepackage{algorithm}
\usepackage[noend]{algpseudocode}

\graphicspath{{./figures/}}

\def\E{\mathbb{E}}

\def\rSigma{\mathrm{\Sigma}} 

\def\hmu{\widehat{\mu}}
\def\hSigma{\widehat{\mathrm{\Sigma}}}

\def\({\left(}
\def\){\right)}
\def\[{\left[}
\def\]{\right]}
\def\nn{\nonumber}

\def\akf{$\alpha$KF~} 

\DeclareMathOperator{\tr}{tr}

\def\wty{\widetilde{y}}

\begin{document}
\title{Nonlinear Kalman Filtering with Divergence Minimization}
\author{San Gultekin and John Paisley\\ Department of Electrical Engineering, Columbia University, New York, NY, USA}
\maketitle

\begin{abstract}
We consider the nonlinear Kalman filtering problem using Kullback-Leibler (KL) and $\alpha$-divergence measures as optimization criteria. Unlike linear Kalman filters, nonlinear Kalman filters do not have closed form Gaussian posteriors because of a lack of conjugacy due to the nonlinearity in the likelihood. In this paper we propose novel algorithms to optimize the forward and reverse forms of the KL divergence, as well as the $\alpha$-divergence which contains these two as limiting cases. Unlike previous approaches, our algorithms do not make approximations to the divergences being optimized, but use Monte Carlo integration techniques to derive unbiased algorithms for direct optimization. We assess performance on radar and sensor tracking, and options pricing problems, showing general improvement over the UKF and EKF, as well as competitive performance with particle filtering.
\end{abstract}


\section{Introduction}
\label{sec:intro}

Modeling and analysis of time-varying signals is one of the most important subfields of signal processing. The problem arises in many different forms, such as communications data sent over a channel, video and audio data, and real-time tracking. A wide variety of algorithms have been developed in the statistics and engineering communities to deal with such dynamic systems. One classic algorithm is the Kalman filter \cite{Kalman_1960}, which performs minimum mean square error estimation of the hidden state of a time-varying linear system. Kalman filter is recursive and online, making it suitable for real-time signal processing applications. Another advantage is its optimality for a large class of state-space models. 

Kalman filtering has been applied extensively in control, communication, and signal processing settings, such as robot motion control and radar target tracking. With the recent explosions in sequential and streaming data, Kalman filters have also become a promising means for approaching machine learning problems, such as natural language processing \cite{Belanger_2015}, collaborative filtering \cite{Gultekin_2014} and topic modeling \cite{Blei_2006}.

An important issue that often arises, requiring modification to basic Kalman filter framework, is nonlinearity. For example, in radar tracking, distance and bearing measurements require a Cartesian-to-polar transformation \cite{LiJilkov3}, whereas dynamic collaborative filtering model contains a bilinear form in two unknown vectors \cite{Koren_2009}. The nonlinear problem has been studied extensively in the literature, resulting in well-known filtering algorithms such as the extended Kalman filter (EKF) \cite{Welch_1995} and unscented Kalman filter (UKF) \cite{Julier_2004}. On the other hand, Monte Carlo methods have been developed \cite{Arulampalam_2002}, which are non-parametric and can represent any probability distribution using a discrete set of points, also referred to as particles. 

While particle filters can approximate arbitrary densities, it may still be important to find the best \textit{parametric} distribution according to a particular objective function. This has been a major goal in Bayesian learning, where the exact posterior distribution is usually intractable and approximated by a known, ``simpler'' distribution. Two established ways to handle this problem are variational inference \cite{Jordan_1999} and expectation-propagation \cite{Minka_2001}, in which the Kullback-Leibler (KL) divergence between the true posterior and the approximating distribution are minimized. Ideas from approximate inference have also been used in the Kalman filtering framework \cite{Beal_2003,Vermaak_2003,Snoussi_2012}. However, a thorough analysis of posterior optimization for nonlinear Kalman filters have not yet been made.

In this paper we fill this gap by presenting three algorithms for nonlinear Kalman filtering based on three respective divergence measures for posterior approximation, each based on a parametric form (in our case, a multivariate Gaussian). These approximations are obtained by algorithms for \textit{approximation-free} divergence minimization. The divergence measures we consider are: 1) the forward KL divergence as used in variational inference; 2) the reverse KL divergence as used in expectation-propagation; and 3) the $\alpha$-divergence, which is a generalized family that contains the former two as special cases. We also show that well-known algorithms such as the EKF and UKF are actually solving \textit{approximations} to KL divergence minimization problems. This further motivates our study to address these shortcomings.\footnote{We emphasize that our methods are still approximate in that the true non-Gaussian posterior will be approximated by a Gaussian. It is approximation-free in that the three algorithms directly optimize the three divergences.} 

The main machinery we use for obtaining these unbiased minimization algorithms is importance sampling. However, the resulting algorithms are all computationally cheaper than particle filtering since 1) no resampling is necessary, and 2) the number of unnecessary samples can be reduced by our proposed adaptive sampling procedure. We show advantages of our algorithms for target tracking and options pricing problems compared with the EKF, UKF and particle filter.

We organize this paper as follows: In Section \ref{sec:background} we define our filtering framework by reviewing the Kalman filter and discussing its non-linear variants. In particular, we discuss parametric approaches, also called assumed density filters, and nonparametric approaches, also called particle filters. In Section \ref{sec:non_kf} we present three divergence minimization problems based on the forward and reverse KL divergence, and $\alpha$-divergence. For each case we propose an algorithm which minimizes the corresponding objective function. Our algorithms are based on Monte Carlo integration techniques. Section \ref{sec:experiments} contains a number of experiments to show how these divergence measures compare with each other and with standard approaches. Finally we conclude in Section 5.

\section{Kalman Filtering}
\label{sec:background}

\subsection{Basic Linear Framework} \label{sec:lin_kf}
The Kalman filter \cite{Kalman_1960} has been developed and motivated as an optimal filter for linear systems. A key property is that this optimality is assured for general state-space models. This has made Kalman filtering widely applicable to a wide range of applications that make linearity assumptions. The Kalman filter can be written compactly at time step $t$ as
\begin{equation}\label{eq:lin_ss}
	x_t = F_t x_{t-1} + w_t , \qquad
	y_t = H_t x_t + v_t ,
\end{equation}
where $w_t$ and $v_t$ are independent zero-mean Gaussian random vectors with covariances $Q_t$ and $R_t$ respectively.\footnote{Kalman's formulation in \cite{Kalman_1960} is optimal for more general noise models, but Gaussian noise is the most common choice, which we also use in this paper.} The latent variable $x_t \in \mathbb{R}^d$ is the unobserved state of the system. The vector $y_t \in \mathbb{R}^p$ constitutes the measurements made by the system. 

The two main tasks of Kalman filtering are prediction and posterior calculation \cite{Welch_1995},
\begin{align} \label{eq:linkf_prob}
	p(x_t | y_{1:t-1}) & = \int p(x_t|x_{t-1})p(x_{t-1}|y_{1:t-1}) dx_{t-1}, \nonumber \\
	p(x_t | y_{1:t}) &\propto p(y_t|x_t)p(x_t | y_{1:t-1}).
\end{align}
When the initial distribution on $p(x)$ is Gaussian all these calculations are in closed form and are Gaussians, which is an attractive feature of the linear Kalman filter. 


\subsection{Nonlinear framework} \label{sec:non_kf}
In many problems the measurements $y_t$ involve nonlinear functions of $x_t$. In this case the Kalman filter becomes nonlinear and the closed-form posterior calculation discussed above no longer applies. The nonlinear process is
\begin{equation}\label{eq:nonlin_ss}
	x_t = F_t x_{t-1} + w_t , \qquad y_t = h(x_t) + v_t ,
\end{equation}
where the noise process is the same as in Eq.\ (\ref{eq:lin_ss}), but $h(\cdot)$ is a nonlinear function of $x_t$.\footnote{We focus on measurement nonlinearity in this paper, assuming the same state space model. The techniques described in this paper can be extended to nonlinearity in the state space as well.} While formally Bayes' rule lets us write
\begin{align}\label{eq:nkfbayes}
	p(x_t|y_{1:t}) = \frac{p(y_t|x_t)p(x_t|y_{1:t-1})}{\int p(y_t|x_t)p(x_t|y_{1:t-1})dx_t},
\end{align}
the normalizing constant is no longer tractable and the distribution $p(x_t|y_{1:t})$ is not known. Although the nonlinearity in $h$ may be required by the problem, a drawback is the loss of fast and exact analytical calculations. In this paper we discuss three related techniques to approximating $p(x_t|y_{1:t})$, but first we review two standard approaches to the problem.

\subsection{Parametric approach: Assumed density filtering} \label{sec:param}
To address the computational problem posed by Eq. \eqref{eq:nkfbayes}, assumed density filters (ADF) project the nonlinear update equation to a tractable distribution. Building on the linear Gaussian state-space model, Gaussian assumed density filtering has found wide applicability \cite{Maybeck_1982,Welch_1995,Julier_2004,Ito_2000,Guo_2006}. The main ingredient here is an assumption of \textit{joint Gaussianity} of the latent and observed variables. This takes the form,
\begin{align} \label{eq:jnt_asmp}
	p(x_t,y_t) \sim N\( \[\begin{array}{c}\mu_x\\\mu_y\end{array}\] , \left[\begin{array}{cc}\rSigma_{xx} & \rSigma_{xy}\\\rSigma_{yx} & \rSigma_{yy}\end{array}\right] \).
\end{align}
(We've suppressed some time indexes and conditioning terms.) Under this joint Gaussian assumption, by standard computations the conditional distribution $p(x_t|y_t)$ is
\begin{eqnarray} \label{eq:jnt_upd}
p(x_t|y_t) &=& N\left(\mu_{x|y} , \rSigma_{x|y}\right),\\[3pt]
 \mu_{x|y} &=& \mu_x + \rSigma_{xy}\rSigma_{yy}^{-1}(y_t-\mu_y),\nonumber\\[3pt]
 \rSigma_{x|y} &=& \rSigma_{xx} - \rSigma_{xy} \rSigma_{yy}^{-1} \rSigma_{yx}.\nonumber
\end{eqnarray}
In this case, the conditional distribution is also the posterior distribution of interest. Using this approximation, Kalman filtering can be carried out. For reference we provide predictive update equations in Appendix \ref{app1}.


There are several methods for making this approximation. We briefly review the two most common here: the extended Kalman filter (EKF) and the unscented Kalman filter (UKF). The EKF approximates $h$ using the linearization 
$$h(x_t) \approx h(x_0) + H(x_0)(x_t-x_0)$$
where $H(x_0)$ is the Jacobian matrix evaluated at the point $x_0$. For example, $x_0$ could be the mean of the prior $p(x_t|y_{1:t-1})$. By plugging this approximation directly into the likelihood of $y_t$, the form of a linear Kalman filter is recovered and a closed form Gaussian posterior can be calculated.

As discussed in \cite{Julier_2004}, the first-order approximation made by the EKF is often poor and performance can suffer as a result. Instead, they propose to estimate the quantities in \eqref{eq:jnt_asmp} with an unscented transform---a numerical quadrature method. The result is the UKF, which has similar computational cost as the EKF and higher accuracy. Based on the calculated Gaussian prior $p(x_t|y_{1:t-1}) = N(x_t|\mu_x,\rSigma_{xx})$, the UKF selects a discrete set of sigma points at which to approximate $\mu_y$, $\rSigma_{yy}$ and $\rSigma_{xy}$.  Let $d_x = \mathrm{dim}(x_t)$ and $N_s = 2\times d_x+1$. These sigma points $x^1,\dots,x^{N_s}$ are
\begin{align}
	x^s &=
	\begin{cases}
		\mu_x &\text{for}~s=0 \\
		\mu_x + [ \sqrt{(d_x+\lambda) \rSigma_{xx}} ]_s &\text{for}~s=1,\ldots,d_x \\
		\mu_x - [ \sqrt{(d_x+\lambda) \rSigma_{xx}} ]_{s-d_x} &\text{for}~s=d_x+1,\ldots,2d_x
	\end{cases}
\end{align}
The vector $[\sqrt{(d_x+\lambda)\rSigma}]_s$ corresponds to the $s$th column of the Cholesky decomposition of the matrix $\sqrt{(d_x+\lambda)\rSigma}$. Positive weights $w^s$ are also defined for each $x^s$. The constant $\lambda$ controls these sigma point locations, as well as the weights (along with additional fixed parameters). These $N_s$ locations are used to empirically approximate all means and covariances in Eq.\ (\ref{eq:jnt_asmp}). Once $y_t$ is measured, the approximation of $p(x_t|y_t)$ can then be calculated using Eq.\ (\ref{eq:jnt_upd}).

There are many extensions to the UKF framework such as cubature Kalman filtering (CKF) \cite{Binjia_2013} and QMC Kalman filtering \cite{Guo_2006}, which use different numerical quadratures to carry out the approximation, but still correspond to the joint Gaussian assumption of Eq.\ (\ref{eq:jnt_asmp}). With that said, however, not all Gaussian ADFs make a joint Gaussianity assumption. For example, methods based on expectation-propagation \cite{Minka_2001} use moment matching (e.g., \cite{Heskes_2002}) to obtain a Gaussian posterior approximation without modifying the joint likelihood distribution. We focus on an EP-like method in Section \ref{sec:reverse}.

\subsection{Nonparametric approach: Particle filtering} \label{sec:nonparam}

We have seen that the main theme of ADF is approximating the posterior with a pre-specified joint probability density; when this joint density is Gaussian then $p(x_t|y_t) \approx N(\mu_t,\rSigma_t)$. On the other hand, nonparametric versions use sampling for posterior approximation without making any density assumptions on the form of this posterior,
\begin{align}\label{eq:par_fil}
	p(x_t|y_{1:t}) ~\approx~ \sum_{s=1}^{N_s} w_t^s \delta_{x_t^s}.	
\end{align} 
The positive weights $w_t^s$ sum to one, and $\delta_{x_t^s}$ is a point mass at the location $x_t^s$. The main approach is to use particle filters, a method base on importance sampling. 
In case of particle filtering using sequential importance resampling (SIR) \cite{Arulampalam_2002}, updating an empirical approximation of $p(x_t|y_{1:t})$ uses a uniform-weighted prior approximation, $p(x_t) \approx \sum_{s=1}^{N_s}\frac{1}{N_s}\delta_{x_{t}^s}$ to calculate the posterior importance weights 
\begin{align}\label{eq:smc}
	w_t^s \propto \frac{1}{N_s}p(y_t|x_t^s).
\end{align}
It then constructs the uniform-weighted prior approximation by sampling $N_s$ times
$$
 x_{t+1}^s \stackrel{iid}{\sim}\, \sum_{s=1}^{N_s} w_t^s N(F_t x_t^s,Q_t),\quad
 p(x_{t+1}) \approx \sum_{s=1}^{N_s} \frac{1}{N_s}\delta_{x_{t+1}^s}
$$
While SIR particle filters can adaptively approximate any posterior density, the double sampling
has computational cost, making these filters considerably slower compared to the above ADF approaches. Another potential issue is the need to propagate particles between time frames, which can be prohibitively expensive in communication-sensitive distributed applications, such as sensor networks \cite{Snoussi_2012}. 


\section{Three divergence minimization approaches}
\label{sec:new_nonkf}

In this section we discuss the three proposed divergence minimization approaches to the nonlinear Kalman filtering problem. These include the two directions of the Kullback-Leibler (KL) divergence as well as the related $\alpha$-divergence that contains both KL divergence measures as limiting cases. In all cases, our goal is to approximate the intractable posterior distribution $p(x_t|y_t)$ with a multivariate Gaussian distribution $q(x_t) = N(\mu_t,\rSigma_t)$, using these three divergences as potential quality measures. We again note that our contribution is to provide three unbiased methods for directly optimizing these divergences without introducing additional approximations. We therefore anticipate an improvement over the standard EKF and UKF approximations. In each of the following three subsections, we present the divergence objective, review its tractability issues, and then present our approach to resolve this issue.

\subsection{Approach 1: Forward KL divergence minimization} \label{sec:forward}

Given two distributions $p(x|y)$ and $q(x)$, the forward KL divergence is defined as
\begin{align} \label{eq:fkl}
	\text{KL}[q\|p] = \int q(x) \ln\frac{q(x)}{p(x|y)} dx .
\end{align}
The KL divergence is always nonnegative, becomes smaller the more $q$ and $p$ overlap, and equals zero if and only if $q=p$. These properties of the KL divergence make it a useful tool for measuring how ``close'' two distributions are. It is not a distance metric however, as $\text{KL}[q\|p] \neq \text{KL}[p\|q]$; we discuss the latter in detail in Section \ref{sec:reverse}. In Bayesian machine learning, minimizing an objective of this form over $q$ is know as variational inference (VI) \cite{Wainwright_2008}. In this case, $p(x|y)$ corresponds to an unknown posterior distribution of the model parameters, and $q$ is its simpler approximation.


For the nonlinear Kalman filtering problem, the posterior is on the latent state vector $x_t$ and so is intractable. Therefore, as is often the case, $\text{KL}[q\|p]$ is not calculable. Variational inference \cite{Bishop_2006,Wainwright_2008} instead uses the identity
\begin{align}
\ln p(y) = \mathcal{L}(q,p(y,x)) + \text{KL}[q(x)\|p(x|y)],
\end{align}
where
\begin{align} \label{eq:elbo}
\mathcal{L}(q,p(y,x)) &= \int q(x) \ln \frac{p(y,x)}{q(x)} dx.
\end{align}
This often is tractable since the joint distribution $p(y,x)$ is defined by the model. Since the marginal $\ln p(y)$ is constant and $\text{KL}\geq 0$, variational inference instead maximizes $\mathcal{L}$ with respect to parameters of $q(x)$ to equivalently minimize KL.

While nonlinear Kalman filters have a simply-defined joint likelihood $p(y_t,x_t|y_{1:t-1})$ at time $t$, a significant problem still arises in calculating $\mathcal{L}$ due to the nonlinear function $h$. That is, if we {define} $q(x_t) = N(\mu_t,\rSigma_t)$, then for the Gaussian generative process of Eq.\ (\ref{eq:nonlin_ss}) we optimize $\mu_t$ and $\rSigma_t$ over the function
\begin{eqnarray}\label{eq:VIobj}
 \mathcal{L} &=& -\frac{1}{2}\mathbb{E}_q[(y_t - h(x_t))^\top R^{-1}(y_t-h(x_t))]\\[2pt]
 &&+ ~\mathbb{E}_q[\ln p(x_t|y_{1:t-1})]-\mathbb{E}_q[\ln q(x_t)] + \text{const.}\nonumber
\end{eqnarray}
The terms in the second line are tractable, but in the first line the nonlinearity of $h(x_t)$ will often result in an integral not having a closed form solution. 

In the variational inference literature, common approaches to fixing this issue typically involve making tractable approximations to $h(x_t)$. For example, one such approximation would be to pick a point $x_0$ and make the first-order Taylor approximation $h(x_t) \approx h(x_0) + H(x_0)(x_t - x_0)$. One then replaces $h(x_t)$ in (\ref{eq:VIobj}) with this approximation and optimizes $q(x_t)$. In fact, in this case the resulting update of $q(x_t)$ is identical to the EKF. This observation implies a correspondence between variational inference and commonly used approximations to the non-linear Kalman filters such as the EKF. We make this formal in the following theorem.\newline

\noindent \textbf{Theorem 1:} Let the joint Gaussian ADF correspond to the class of filters which make the joint distribution assumption in \eqref{eq:jnt_asmp}. Then, all filters in this class optimize an approximate form of the variational lower bound in \eqref{eq:VIobj}.\newline

We present a complete proof in Appendix \ref{app2}. Theorem 1 is general in that it contains the most successfully-applied ADFs such as the EKF and UKF, among others. For the special case of EKF, the nature of this approximation is more specific.\newline

\noindent \textbf{Corollary 2:} The EKF corresponds to optimizing the objective \eqref{eq:VIobj} using a first order Taylor approximation of $h$.\newline

Please see Appendix \ref{app3} for a proof. Consequently, the existing algorithms \textit{modify} $\mathcal{L}$ and the optimization of this approximation to $\mathcal{L}$ over the parameters of $q(x_t)$ is no longer guaranteed to minimize $\text{KL}[q\|p]$. Instead, in this paper we are motivated to fill in this gap and find ways to \textit{directly} optimize objectives such as \eqref{eq:VIobj}, and thus minimize divergence measures between $q$ and the intractable posterior $p(x_t|y_t)$.
We next devise a method for $\text{KL}[q\|p]$.

Recently Paisley, et al. \cite{Paisley_2012} proposed a stochastic method for sampling \textit{unbiased} gradients of $\mathcal{L}$, allowing for approximation-free minimization of the forward KL divergence using stochastic gradient descent. We derive this technique for the nonlinear Kalman filter, which will allow for approximate posterior inference having smaller KL divergence than the EKF and UKF.
Using simpler notation, we seek to maximize an objective of the form,
\begin{eqnarray}
 \mathcal{L} &=& \mathbb{E}_q[f(x_t)] + \mathbb{E}_q[\ln p(x_t)] - \mathbb{E}_q[\ln q(x_t)]~~ \label{eq:elbo_kf} \\[5pt]
 f(x_t) &=& -\frac{1}{2}(y_t - h(x_t))^\top R^{-1}(y_t - h(x_t))\label{eq.f}
\end{eqnarray}
over the parameters of $q(x_t) = N(x_t|\mu_t,\rSigma_t)$, and thereby minimize $\text{KL}[q\|p]$. This can be done by gradient ascent. However, since $\mathbb{E}_q[f(x_t)]$ does not have a closed form solution, $\nabla\mathcal{L}$ can not be evaluated analytically. The proposed solution in \cite{Paisley_2012} is to instead step in the direction of an unbiased stochastic gradient. To this end, the observation is made that
\begin{eqnarray}
 \nabla \mathcal{L} = \mathbb{E}_q[f(x_t)\nabla\ln q(x_t)] + \nabla\mathbb{E}_q\Big[\ln \frac{p(x_t)}{q(x_t)}\Big],
\end{eqnarray}
where the identity $\nabla q(x_t) = q(x_t)\nabla\ln q(x_t)$ is used. While the second gradient can be calculated analytically with respect to either $\mu_t$ or $\rSigma_t$, the first gradient can be sampled using Monte Carlo integration,
\begin{equation}
 \mathbb{E}_q[f(x_t)\nabla\ln q(x_t)] \approx \frac{1}{S}\sum_{s=1}^S f(x^s_t)\nabla\ln q(x^s_t),\quad x_t^s \stackrel{iid}{\sim} q(x_t).
\end{equation}
A second observation is made by \cite{Paisley_2012} that the variance of these samples may be so large that $S$ needs to be set to too large a number to make this approximation computationally feasible.
For this reason employing variance reduction methods is crucial. Paisley, et al.\ \cite{Paisley_2012} propose introducing a \textit{control variate} $g(x_t)$ that is highly correlated with $f(x_t)$, but has an analytic expectation $\mathbb{E}_q[g(x_t)]$. The gradient of $\mathcal{L}$ with a control variate is equal to
\begin{eqnarray}
 \nabla \mathcal{L} &=& \mathbb{E}_q[(f(x_t)-\lambda g(x_t))\nabla\ln q(x_t)] + \lambda\nabla\mathbb{E}_q[g(x_t)]\nonumber\\[5pt]
 &&+ ~\nabla\mathbb{E}_q[\ln p(x_t)] - \nabla\mathbb{E}_q[\ln q(x_t)].
\end{eqnarray}
Though this leaves the gradient unchanged, MC sampling of the first term has much smaller variance when $|corr(f,g)|$ is large (calculated using $q(x_t)$). The parameter $\lambda\in\mathbb{R}$ is set to minimize the variance.\footnote{As shown in \cite{Paisley_2012}, when $\lambda \equiv cov(f,g)/var(g)$ (approximated by sampling) the variance reduction ratio is $var(f-\lambda g)/var(f) = 1-corr(f,g)^2$.} Intuitively, this can be seen by noting that if $f(x_t^s) \approx \lambda g(x_t^s)$ at the sampled values $x_t^s$, then $|f(x_t^s) - \lambda g(x_t^s)| \ll |f(x_t^s)|$. In this case, the {analytic} gradient $\lambda\nabla\mathbb{E}_q[g(x_t)]$ gives an initial approximation of $\mathbb{E}_q[f(x_t)\nabla \ln q(x_t)]$, which is then corrected to be made unbiased by the MC-sampled $\mathbb{E}_q[(f(x_t)-\lambda g(x_t))\nabla\ln q(x_t)]$. Since $g(x_t)$ is a good approximation of $f(x_t)$ in the region of high probability defined by $q(x_t)$, the analytic approximation captures most information, but is refined by 
the MC-sampled gradient to make the method approximation-free.

\begin{algorithm}[t]
\caption{SKF: stochastic search Kalman filter}\label{alg1}
\begin{algorithmic}[1]
\State \textbf{Input:} Posterior $q(x_{t-1})$, sample size $S$, and iterations $I$.
\State Calculate prior $p(x_t) = N(\mu_t,\rSigma_t)$.
\For{$i=1,\dots,I$}
\State Sample $x_t^s \sim_{iid} q(x_t)$ for $s=1,\dots,S$. 
\State Compute $\nabla_{\widehat{\mu}_t}\mathcal{L}$ and $\nabla_{\widehat{\rSigma}_t}\mathcal{L}$ as in \eqref{eq:sgrd1} and \eqref{eq:sgrd2}.
\State Update 
\begin{eqnarray}
 \widehat{\mu}_t^{(i+1)} &=& \widehat{\mu}_t^{(i)} + \rho_i ~ [ C^{(i)} ~ \nabla_{\widehat{\mu}_t} \mathcal{L} ]\nn \\
 \widehat{\rSigma}_t^{(i+1)} &=& \widehat{\rSigma}_t^{(i)} + \rho_i ~ [ C^{(i)} ~ \nabla_{\widehat{\rSigma}_t} \mathcal{L} ~ {C^{(i)}} ] \nn
 \end{eqnarray}  
\EndFor 
\State \textbf{Return} $q(x_t) = \mathcal{N}(\widehat{\mu}_t,\widehat{\rSigma}_t)$
\end{algorithmic}
\end{algorithm}

The requirements on $g(x_t)$ to be a good control variate for $f(x_t)$ are that: 1) it is an approximation of $f(x_t)$, and 2) the expectation  $\mathbb{E}_q[g(x_t)]$ is solvable. There are many possible control variates for the function $(y-h(x))^\top R^{-1} (y-h(x))$. However, building on the EKF framework we propose setting 
\begin{eqnarray}
 g(x_t) &=& -\frac{1}{2}(y_t - \widetilde{h}(\mu_t,x_t))^\top R^{-1}(y_t - \widetilde{h}(\mu_t,x_t))\nn\\[5pt]
 \widetilde{h}(\mu_t,x_t) &=& h(\mu_t) + H(\mu_t)(x_t-\mu_t)
\end{eqnarray}
We let $\mu_t$ be the current value of the mean of $q(x_t)$ at a given iteration of time $t$. If we define $\widetilde{y}_t = y_t - h(\mu_t) + H(\mu_t)\mu_t$, then equivalently we can write
\begin{equation}
 g(x_t) = -\frac{1}{2}(\widetilde{y}_t - H(\mu_t)x_t)^\top R^{-1} (\widetilde{y}_t - H(\mu_t)x_t).
\end{equation}
The expectation is now in closed form. While a better approximation may have greater variance reduction for a fixed number of MC-samples, we emphasize this would not make the algorithm more ``correct.'' Where the EKF simply \textit{replaces} $f(x)$ with $g(x)$, our stochastic gradient approach then \textit{corrects} the error of this approximation.

Next we derive the unbiased gradients. In this case, these gradients are $\nabla_{\hmu_t}\mathcal{L}$ and $\nabla_{\hSigma_t}\mathcal{L}$. We note that, if we were using the EKF framework, by replacing $f(x_t)$ with $g(x_t)$ in Eq.\ (\ref{eq.f}), the roots of these gradients could be solved and the EKF solutions for $\hmu_t$ and $\hSigma_t$ would result. However, since we have the additional stochastic gradient term, we must perform gradient ascent. The final expressions for the unbiased gradients using samples $x_t^s \sim_{iid} q(x_t)$ are: 
\begin{align}
\nabla_{\widehat{\mu}_t} \mathcal{L} &= \frac{1}{S} \sum_{s=1}^{S} \big[ f(x^s) - g(x^s) \big] \big[ \widehat{\rSigma}_t^{-1} x^s - \widehat{\rSigma}_t^{-1} \widehat{\mu}_t \big] \label{eq:sgrd1} \\
&\qquad + \rSigma_t^{-1} (\mu_t -  \widehat{\mu}_t) + H(x_0)^\top R^{-1}( \wty_t-H(x_0) \widehat{\mu}_t)\nn,  \\
\nabla_{\widehat{\rSigma}_t} \mathcal{L} &= \frac{1}{S} \sum_{s=1}^{S} \left[ f(x^s) - g(x^s) \right] \times \nn \\
&\quad \frac{1}{2}\big[ \widehat{\rSigma}_t^{-1} (x^s-\widehat{\mu}_t) (x^s-\widehat{\mu}_t)^\top \widehat{\rSigma}_t^{-1} - \widehat{\rSigma}_t^{-1} \big] \nonumber \\
&\quad + \frac{1}{2}( \widehat{\rSigma}_t^{-1} - \rSigma_t^{-1}) - \frac{1}{2} H(x_0)^\top R^{-1} H(x_0) ~. \label{eq:sgrd2}
\end{align}

Comparing the gradients \eqref{eq:sgrd1}-\eqref{eq:sgrd2} we see that stochastic search acts as a correction step to the EKF updates resulting from the last lines alone. On the other hand, unlike EKF, we cannot simply solve for $\widehat{\mu}_t$ and $\widehat{\rSigma}_t$ by settings the gradients equal to zero, so we use gradient ascent. Without proper scaling we can easily have a numerically unstable algorithm and the covariance matrix can lose its positive definiteness. To fix this we pre-condition the gradients with a symmetric positive definite matrix $C$ and perform the following updates
\begin{align}
\widehat{\mu}_t^{(i+1)} &= \widehat{\mu}_t^{(i)} + \rho^{(i)} ~ [ C^{(i)} ~ \nabla_{\widehat{\mu}_t} \mathcal{L} ] \label{eq:ngrd1}, \\
\widehat{\rSigma}_t^{(i+1)} &= \widehat{\rSigma}^{(i)} + \rho^{(i)} ~ [ C^{(i)} ~ \nabla_{\widehat{\rSigma}_t} \mathcal{L} ~ {C^{(i)}}] ~. \label{eq:ngrd2}	
\end{align}

We note the difference between index $t$ and $i$. The first is the time frame we are currently processing, while the second is the iteration number at time $t$ since we are using a gradient optimization method. For the conditioning matrix we choose $C^{(i)} = [\rSigma_t^{(i)}]^{-1}$. Using this setting, we get approximate natural gradients \cite{Amari_1998} for $\widehat{\mu}$ and $\widehat{\rSigma}$.
When the step size satisfies the Robbins-Monro conditions, $\sum_{i=1}^{\infty} \rho^{(i)} = \infty$ and $\sum_{i=1}^{\infty} [\rho^{(i)}]^2 < \infty$, the gradients in \eqref{eq:ngrd1}-\eqref{eq:ngrd2} converge to a fixed point of the exact variational lower bound. In practice we can, for example, choose $\rho^{(i)} = (w+i)^{-\eta}$ with $\eta \in (0.5,1]$ and $w \geq 0$. In simulations we observed that, when natural gradients are employed a generic schedule for step sizes can be used, and no 
further hand-tuning is necessary. We refer to this algorithm as \textit{stochastic search Kalman filtering} (SKF) and summarize it in Algorithm 1 for a single time step.

\subsection{Approach 2: Reverse KL divergence minimization}
\label{sec:reverse}

As mentioned in Section \ref{sec:forward}, KL divergence is not a distance measure since it is not symmetric. The complement of the forward KL divergence defined in \eqref{eq:fkl} is the reverse KL divergence:
\begin{align} \label{eq:rkl}
	\text{KL}[p\|q] = \int_{\mathcal X} p(x) \log\frac{p(x)}{q(x)} dx .
\end{align}
We can see that \eqref{eq:rkl} offers an alternative measure of how similar two probability distributions are; therefore we can use it to approximate an intractable posterior distribution.

\begin{algorithm}[t]
\caption{MKF: moment matching Kalman filter}\label{alg2}
\begin{algorithmic}[1]
\State \textbf{Input:} Posterior $q(x_{t-1})$, sample size $S$, proposal dist.\ $\pi_t$.
\State Calculate prior $p(x_t) = N(\mu_t,\rSigma_t)$.
\State Sample $x^s \sim_{iid} \pi_t(x_t)$ for $s=1,\dots,S$. 
\State Calculate $w^s = p(y_t|x^s)p(x^s)/\pi_t(x^s)$, $W = \sum_{s=1}^S w^s$.
\State Approximate the moments of $p(x_t|y_t)$	as
\begin{eqnarray}
 \hmu_t &=& \textstyle\frac{1}{W}\sum_{s=1}^S w^s x^s \nn\\
 \hSigma_t &=& \textstyle\frac{1}{W}\sum_{s=1}^S w^s(x^s-\hmu_t)(x^s-\hmu_t)^\top\nn
\end{eqnarray}

\State \textbf{Return} $q(x_t) = \mathcal{N}(\hmu_t,\hSigma_t)$
\end{algorithmic}
\end{algorithm}

Note that for either objective function, \eqref{eq:fkl} or \eqref{eq:rkl}, the optimal solution will be $q(x)=p(x|y)$. However, since the approximating distribution is typically different from the exact posterior distribution, the two optimization problems will give different solutions in practice. In particular, reverse KL divergence has shown to be a better fit for unimodal approximations, while forward KL works better in multimodal case \cite{Bishop_2006}. Consequently, we can expect that optimizing the reverse KL will be a better choice for the nonlinear Kalman filtering problem (this is supported by our experiments). In Section \ref{sec:forward}, finding a fixed point of the forward KL problem required an iterative scheme for maximizing the variational objective function. The fixed point of the reverse KL has a more interpretable form, as we will show.

To this end, we first note that an exponential family distribution has the form 
$$q(x) = h(x) \exp \{ \eta^\top s(x) - \log A(\eta) \},$$
where $\eta$ is the natural parameter and $s(x)$ is the sufficient statistic. Therefore inference in exponential families correspond to determining $\eta$. Substituting this parametrized form in \eqref{eq:rkl} and setting the derivative with respect to the natural parameter equal to zero, one can show that
\begin{equation}
	0 ~=~ \nabla_{\eta} \text{KL}[p\|q] ~=~ \E_q[s(x)] - \E_p[s(x)], \nn	
\end{equation}
which follows from the exponential family identity $\nabla_{\eta} \log A(\eta) = \E_q[s(x)]$. Therefore the fixed points of the objective are given by
\begin{equation} \label{eq:mm}
	\E_q [s(x)] = \E_p [s(x)], 
\end{equation}
This moment matching is well-known in statistics, machine learning, and elsewhere \cite{Bishop_2006}. In machine learning it appears prominently in expectation propagation \cite{Minka_2001,Lobato_2016}.

A common choice for the approximating exponential family distribution is again Gaussian because it is the maximum entropy distribution for the given first and second order moments \cite{Wainwright_2008}. Since a Gaussian is completely specified by its mean and covariance, when the approximating distribution $q(x)$ is selected to be Gaussian, the optimal solution is simply found by matching its mean and covariance to that of $p(x|y)$.

Therefore, in the context of exponential families the task of finding the optimal distribution for the reverse KL reduces to the task of matching moments. However, there is still a difficulty in the need to compute the moments of an unknown posterior distribution. Fortunately, Monte Carlo methods prove useful here as well. Let $I(f) = \E_q(f(x))$ be the expectation we wish to compute. For example, choosing $f(x) = x$ and $f(x) = xx^\top - \E[x]\E[x]^\top$ gives the mean and covariance respectively. This expectation can be approximated as 
\begin{align}
	\E_q[f(x)] &= \int f(x) \frac{p(x|y)}{\pi(x)} \pi(x) dx , \nn \\
	&= \int f(x) \frac{[p(y|x)p(x)]/\pi(x)}{\int p(y|x')p(x')dx'} \pi(x) dx , \nn \\
	& \approx \sum_{s=1}^S f(x^s) \frac{[p(y|x^s)p(x^s)]/\pi(x^s)}{\sum_j [p(y|x^j)p(x^j)]/\pi(x^j)}. \label{eq:imp_mm}
\end{align}
We will define $$w^s = \frac{p(y|x^s)p(x^s)}{\pi(x^s)},\quad W = \sum_s \frac{p(y|x^s)p(x^s)}{\pi(x^s)},$$
and so $\E_q[f(x)] \approx  \frac{1}{W}\sum_{s=1}^S f(x^s) w^s$. This is related to importance sampling, with the added normalizer $W$. As we can see from \eqref{eq:imp_mm} this procedure is biased as it is a ratio of two approximations, yet it converges to the true expectation $\E_q[f(x)]$ almost surely. Therefore, we have an asymptotically unbiased divergence minimization procedure.
We call this the \textit{moment matching Kalman filter} (MKF) and summarize it in Algorithm 2 for Gaussian distributions, as Gaussian approximations are our focus in this paper. We observe that a major difference between the MKF and SKF of the previous section is that the MKF only needs to sample once to obtain the moment estimates for a time step. Therefore, the MKF is not an iterative algorithm and is much faster. Also, the MKF is slightly faster than particle filtering as it eliminates the need for resampling.

\subsection{Approach 3: $\alpha$-divergence minimization}
\label{sec:alpha}

In Sections \ref{sec:forward} and \ref{sec:reverse} we showed how nonlinear Kalman filtering can be performed by minimizing the forward and reverse KL divergence. 
A further generalization is possible by considering the $\alpha$-divergence, which contains both KL divergences as a special case. Following \cite{Lobato_2016} we define the $\alpha$-divergence to be
\begin{align} \label{eq:adiv}
	D_{\alpha}[p\|q] = \frac{1}{\alpha(1-\alpha)} \( 1 - \int p(x)^\alpha q(x)^{1-\alpha} \) ,
\end{align}
where the parameter $\alpha$ can take any value in $(-\infty,\infty)$. Some special cases are
\begin{align} \label{eq:alp_cases}
	&\underset{\alpha \rightarrow 0}{\lim} D_{\alpha}[p\|q] = \text{KL}[q\|p] ~,~ \underset{\alpha \rightarrow 1}{\lim} D_{\alpha}[p\|q] = \text{KL}[p\|q] ~, \nn \\
	&D_{\frac{1}{2}}[p\|q] = 2 \int (\sqrt{p(x)} - \sqrt{q(x)})^2 dx = 4 \text{Hel}^2[p\|q] ~,
\end{align}
where $\text{Hel}[p\|q]$ is the Hellinger distance. We see that when $\alpha = 1/2$ we get a valid distance metric. Similar as before, we now seek a $q$-distribution which approximates $p(x|y)$, where approximation quality is now measured by the $\alpha$-divergence. 

\begin{algorithm}[t]
\caption{\akf: $\alpha$-divergence Kalman filter}\label{alg3}
\begin{algorithmic}[1]
\State \textbf{Input:} Posterior $q(x_{t-1})$, sample size $S$, and proposal $\pi_t$.
\State Calculate prior $p(x_t) = N(\mu_t,\rSigma_t)$.
\State Sample $x^s \sim_{iid} \pi_t(x_t)$ for $s=1,\dots,S$. 
\State Calculate $w^s = \frac{[p(y_t|x^s)p(x^s)]^\alpha q(x^s)^{1-\alpha}}{\pi_t(x^s)},~ W = \sum_{s=1}^S w^s$
\State Approximate the moments of $\widetilde{p}(x_t)$
 \begin{eqnarray*}
  \hmu_t &=& \textstyle\frac{1}{W}\sum_{s=1}^S w^s x^s \\
  \hSigma_t &=& \textstyle\frac{1}{W}\sum_{s=1}^S w^s(x^s-\hmu_t)(x^s-\hmu_t)^\top
 \end{eqnarray*}
\State \textbf{Return} $q(x_t) = \mathcal{N}(\hmu_t,\hSigma_t)$
\end{algorithmic}
\end{algorithm}

Again assuming that the approximating distribution is in the exponential family, $q(x) = h(x) \exp \{ \eta^\top s(x) - \log A(\eta) \}$. The gradient of the $\alpha$-divergence shows that
\begin{align} \label{eq:alp_der}
0 & \,=\, \nabla_\eta D_{\alpha}[p\|q]
\,=\, \frac{1-\alpha}{Z_{\widetilde{p}}} \int \widetilde{p}(x) \big[ s(x) - \E_q[s(x)] \big] \nn \\
&=\, \E_{\widetilde{p}}[s(x)] - \E_q[s(x)] 
\end{align}
Note that we defined a new probability distribution $\widetilde{p}(x) = p(x)^\alpha q(x)^{1-\alpha} / Z_{\widetilde{p}}$ where the denominator term is the cumulant function. This leads to a \emph{generalized} moment matching condition,
\begin{equation} \label{eq:gmm}
	\E_q [s(x)] = \E_{\widetilde{p}} [s(x)] .
\end{equation}
This problem is more complicated than the reverse KL because the left hand side also depends on the $q$-distribution. 
The $\alpha$-divergence generalizes a number of known divergence metrics. In context of EP, it is possible to obtain a generalization which is called Power EP \cite{Minka_2004}. More recently, \cite{Lobato_2016} used a similar black-box optimization, where they showed that by varying the value of $\alpha$ the algorithm varies between variational inference and expectation propagation. It turns out that, for many practical problems, using a fractional value of $\alpha$ can give better performance than the limiting cases $\alpha=0$ or $\alpha=1$. This motivates our following $\alpha$-divergence minimization scheme.

A similar importance sampling methodology can be used for this optimization as for the reverse KL divergence. Using similar notation, we can write
\begin{align}
	\E_{\widetilde{p}}[f(x)] &= \int f(x) \frac{p(x|y)^\alpha q(x)^{1-\alpha}}{\pi(x)} \pi(x) dx ~, \label{eq:imp_gmm} \\
	& \approx \sum_{s=1}^S f(x^s) \frac{[p(y|x^s)p(x^s)]^\alpha q(x^s)^{1-\alpha}/\pi(x^s)}{\sum_j [p(y|x^j)^\alpha p(x^j)^\alpha q(x^j)^{1-\alpha}]/\pi(x^j)},\nn
\end{align}
where $x^s \sim_{iid} q(x)$. Again we define 
$$w^s = [p(y|x^s)p(x^s)]^\alpha q(x^s)^{1-\alpha}/\pi(x^s), ~~ \textstyle W=\sum_s w^s.$$ We see that the procedure in \eqref{eq:imp_mm} is a special case of this when we set $\alpha=1$. However, there is a significant difference in that the moment matching of \eqref{eq:mm} can be done in one iteration since it only depends on $p$. In \eqref{eq:imp_gmm} the $q$ distribution appears on both sides of the equality. 
%
%
%
%
%
This is similar to of EP and Power-EP algorithms, where multiple iterations can be run to update $q$. Upon convergence we know that the solution is a fixed point of \eqref{eq:adiv}, but convergence of the procedure is not guaranteed and multiple iterations might degrade the performance. In our experiments we will only iterate once to avoid possible diverging and also to keep the cost of the algorithm the same as that of MKF in the previous section. We call this algorithm $\alpha$-divergence Kalman filter ({\akf}) and summarize it in Algorithm 3. We note that the only difference between \akf and MKF is in step 4.

We can get a better understanding of employing $\alpha$-divergence by analyzing the weight coefficients. In particular, lets assume that we choose our proposal distribution as the prior, i.e. $\pi(x)=p(x)$. Then, the MKF weights become $w^s \propto p(y|x^s)$ in the Kalman filter. The $\alpha$KF weights, on the other hand are $w^s \propto p(y|x^s)^\alpha$; therefore, the likelihood term is scaled by alpha and as $\alpha \rightarrow 0$ all the particles generated will have equal contribution. For very low values of $\alpha$ this will discard all the information, which is clearly unwanted, but for intermediate values this can alleviate the effects of sharply fluctuating likelihood factors. As we will show in our experiments, when the measurement noise is strong, choosing an intermediary $\alpha$ value provides robustness.

\subsection{Adaptive Sampling} \label{sec:adasamp}

\begin{figure}[t]
	\centering
	\includegraphics[scale=.7]{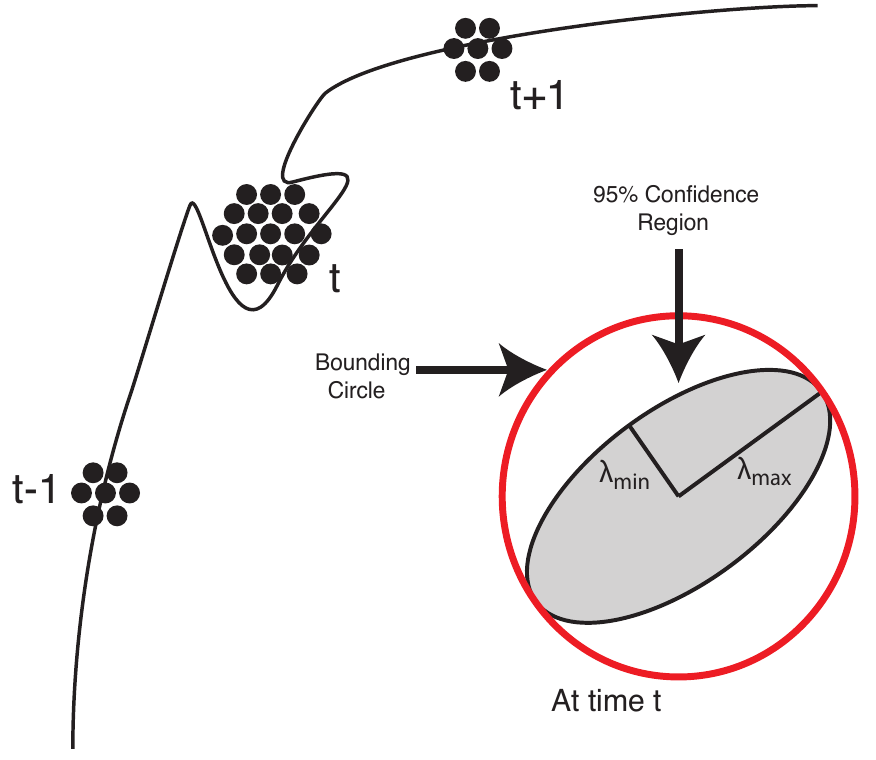}
	\caption{Illustration of adaptive sampling. Due to unexpected changes in a target trajectory, more samples may be needed at a given time point. Also shown is the bounding circle for a confidence ellipsoid in two dimensions.}
	\label{fig:adasamp}
\end{figure}

The main parameter in the implementation of sampled filters such as particle filters and the three filters proposed here is the number of particles $S$ that will be used. Hence it is desirable to have a method of estimating the minimum number of samples necessary for a given degree of accuracy. Then, for each round of filtering we can use this computed sample size to reduce the computation as much as possible, but still be able to increase the sample size when necessary. In Figure \ref{fig:adasamp} we illustrate the problem of tracking a moving target. At time $t$ this target makes an abrupt maneuver where we need more particles for accurate tracking, but we can reduce the size afterwards. 

Since the approximate Gaussian posterior distribution has a parametric form, we are able to use an adaptive sampling method for the MKF and \akf.\footnote{For the SKF the per-iteration sample size is much smaller, so there is less benefit in using this technique.} To determine the appropriate number of samples, we measure the uncertainty of our mean approximation for $q(x_t) = \mathcal{N}(\mu_t,\rSigma_t)$, where $\mu_t = \sum_{s=1}^S w^s x^s/W$. For importance sampling, the variance of this estimator is approximately
\begin{align}
	\mathbb{V}(\mu_t) \,\approx\, \sum_{s=1}^S \[ \frac{w^s}{W} \]^2 (x^s - \mu_t) (x^s - \mu_t)^\top,
\end{align} 
Here, if $S$ is large enough the estimator can be approximated as normal by the central limit theorem \cite{Andrieu_2003}. We use this to compute the radius of a 95\% confidence region. Without loss of generality assume that the estimator is zero mean, which is justified by the asymptotic unbiasedness of the unnormalized importance sampling procedure. We denote the estimator by $\widehat{X} \sim N(0,\mathbb{V}(\mu_t))$. Then we have
\begin{align} \label{eq:con_reg}
	P( \widehat{X}^\top \mathbb{V}(\mu_t)^{-1} \widehat{X} \leq \chi_d^2(p) ) = p ,
\end{align}
where $\chi_d^2(p)$ is the quantile function of chi-squared distribution with d degrees of freedom (which equals the state-space dimension here), and $p$ is the probability value (for 95\% confidence intervals this is set to $p=0.95$). $\chi_d^2$ indicates the chi-squared distribution. The region described by \eqref{eq:con_reg} is a hyper-ellipsoid, so the maximum possible radius will correspond to the major axis, which is given by 
\begin{align}\label{eq:rmax}
	r_\text{max} = \sqrt{ \lambda_\text{max}(\mathbb{V}(\mu_t)) \times \chi_d^2(0.95) } .
\end{align}
Note that this is a conservative estimate, as the hypersphere with radius $r_\text{max}$ will typically be much larger than the hyper-ellipsoid. An illustration of the \emph{bounding circle} for 2D multivariate normal distribution is given in Figure \ref{fig:adasamp}.

Now assume that using a small sample set $S_\text{base}$ we wish to estimate the minimum number of samples $S_\text{min}$ required to achieve a certain $r_\text{max}$. 
We have the relation and result that
\begin{equation} \label{eq:rmax_ratio}
	\frac{r_1}{r_2} \propto \sqrt{ \frac{S_2}{S_1} },\qquad S_\text{min} = S_\text{base} \times \[\frac{r_\text{base}}{r_\text{max}}\]^2 .
\end{equation}

As expected, the smaller radius we desire, the larger sample size we need. We note in passing that $\mathbb{V}(\mu_t)$ is our confidence in estimating the mean of the true posterior, and not the ground truth. The accuracy of estimating the latter is dictated by the measurement noise, and cannot be made arbitrarily small by increasing the sample size.

\section{Numerical Results} \label{sec:experiments}
We experiment with all three proposed nonlinear Kalman filter in algorithms, as well as the EKF, UKF and particle filter, on radar and sensor tracking problems, as well as an options pricing problem.

\subsection{Target Tracking}

The first problem we consider is target tracking. This problem arises in various settings, but here we consider two established cases: radar and sensor networks. The radar tracking problem has been a primary application area for nonlinear Kalman filtering. The target is typically far away from the radar, for example an airplane. Wireless sensor networks are another emerging area where nonlinear filtering is useful. Driven by the advances in wireless networking, computation and micro-electro-mechanical systems (MEMS), small inexpensive sensors can be deployed in a variety of environments for many applications \cite{Chong_2003,Tubaishat_2003}. 

For both problems the state-space will have the form
\begin{eqnarray}
x_t &= ~F_t x_{t-1} + w_t ,&~ w_t \sim N(0,Q_t), \nn \\
y_t &\hspace{-9pt}= ~h(x_t) + v_t ,&~ v_t \,\sim N(0,R_t) \label{eq:radar_model}.
\end{eqnarray}

Here, $F_t$ and $Q_t$ model the dynamics of target motion and are usually time-varying. On the other hand, $h(\cdot)$ specifies the equipment that performs the measurements, and the environment and equipment based inaccuracies are represented by $R_t$. In the radar setting, when the target is far away and the angle measurement noise is strong enough, the problem can become highly nonlinear. For sensor networks, the nonlinearity is caused by the small number of active sensors (due to energy constraints) with large measurement noise (due to the attenuation in received signal) \cite{Boukerche_2007}.  While the value of $R_t$ can be determined to some extent through device calibration, it is more challenging to do this for $Q_t$ \cite{LiJilkov5}.

Our experiments are based on synthetic data using a constant velocity model in $\mathbb{R}^2$ which corresponds to the state vector vector $x_t = [x_1, ~ \dot{x_1}, ~ x_2, ~ \dot{x_2}]^{\top}$; the second and fourth entries correspond to the velocity of the target in each dimension. Following \cite{LiJilkov1}, we set the parameters for the state variable equation to 
\begin{align}
\label{eq:model_CV}
F_t = \,\,
\begin{bmatrix}
F_2 & 0 \\ 0 & F_2
\end{bmatrix}, &~
~~F_2 = 
\begin{bmatrix}
1 & \mathrm{\Delta} t \\ 0 & 1
\end{bmatrix}, \\
Q_t =
\begin{bmatrix}
Q_2 & 0 \\ 0 & Q_2
\end{bmatrix}, &~
~~Q_2 = \sigma_{CV}
\begin{bmatrix}
\mathrm{\Delta} t^4/4 & \mathrm{\Delta} t^3/2 \\ \mathrm{\Delta} t^3/2 & \mathrm{\Delta} t^2
\end{bmatrix}.
\end{align} 

The radar measures the distance and bearing of the target via the nonlinear function $h(\cdot)$ of the target location, $$h(x_t) = \big[ \sqrt{x_t(1)^2 + x_t(3)^2},~~ \tan^{-1} [x_t(3)/x_t(1)] \big]^\top, $$ i.e. the Cartesian-to-polar transformation \cite{LiJilkov3}. For the sensor networks, we will consider a scenario which uses range-only measurements from multiple sensors. This yields the model in \eqref{eq:radar_model} where $h(\cdot)$ is the measurement function such that the $i$-th dimension (i.e. measurement of sensor $s_i$) is given by $$[h(x_t)]_i = \sqrt{[x_t(1)-s_i(1)]^2 + [x_t(2)-s_i(2)]^2},$$ and the length of $h(x_t)$ will be the number of activated sensors at time $t$.

We consider two types of problems: tracking with uncertain parameters and tracking with known parameters. For the case of uncertain parameters, we set the radar and sensor simulation settings as follows. First, we note that for both simulations we assume a constant measurement rate, and so set $\Delta t = 1$. For radar we sweep the process noise values in \eqref{eq:model_CV} as $\sigma_{CV} \in \{10^{-3}, ~2 \times 10^{-3}, ~\ldots~ , ~ 10^{-2} \}$. We generate 20 data sets for each value of $\sigma_{CV}$, yielding a total of 200 experiments. For the measurement noise we use a diagonal $R$ with entries $\sigma_r^2=10^{-1}$ and $\sigma_{\theta}^2=10^{-2}$ which dictates the noise of distance and bearing measurements respectively. The initial state is selected as $x_0 = [1000,~ 10,~ 1000,~ 10]^\top$; this distance from origin and angle noise variance results in a severely nonlinear model, making filtering quite challenging. For sensor network simulations, we use the same constant-velocity model of 
\eqref{eq:model_CV} with $\sigma_{CV}=10^{-2}$. We deploy 200 sensors and at each time there are exactly 3 distinct ones responsible for range measurements. The measurement covariance matrix is $R = \sigma_R^2 I$ where we set $\sigma_R = 20$. We select the initial state as $x_0 = [1000,~ 1,~ 1000,~ 1]^\top$. With this, once again, we obtain a highly nonlinear system, albeit less severe than the radar case. We also consider the case where the generating parameters are known to the filter. In this case, we assign the performance of the filter as a function of process and measurement noise covariances. For this one, we sweep $\sigma_\text{CV} \in \{0.001, 0.005, 0.01, 0.05, 0.1\}$ and $\sigma_r \in \{10,15,20,25,30\}$. We report the results for the sensor network case.

We implemented EKF, UKF, sampling-importance-resampling particle filter (PF), and our proposed SKF, MKF, and \akf for $\alpha=0.5$. For SKF we use $500$ particles/iteration, whereas we consider $10^4$ particles for PF, MKF and \akf. When there is parameter uncertainty, the exact value of $Q$ is not known to the filter, therefore we consider a scaled isotropic covariance of form $\sigma_Q^2 I$.

\begin{table}[t]
	\caption{Radar tracking problem: Mean Square Error (MSE) of various filtering schemes as a function of process noise parameter $\sigma_Q$. The boldfaces show the best performers for small/large particle sizes.}
	\label{table_radar}
	\centering
	 \def\arraystretch{1.25}
	\begin{tabular}{ |c|c|c|c|c|c| } 
		\hline
		& \multicolumn{5}{|c|}{$\sigma_Q$} \\ \cline{2-6}
		&$10^{-2}$ & $5 \times 10^{-2}$ & $10^{-1}$ & $5 \times 10^{-1}$ & $1$ \\
		\hline
		SKF           & 41.4100 & 34.6611 & 29.9952 & 42.1360 & 38.0507 \\
		MKF & 31.3088 & 27.6861 & 29.0376 & 35.2422 & 39.2536 \\ 
		$\alpha$KF    & 30.8783 & 27.9475 & \textbf{27.4130} & 31.0271 & 34.9420 \\
		\hline
		PF            & 28.5429 & 32.3768 & 35.1842 & 44.3704 & 48.9767 \\ 
		\hline
		EKF           & 33.8611 & 35.8086 & 37.7808 & 42.6595 & 45.9788 \\
		UKF           & 31.7528 & 31.8616 & 33.7625 & 41.1282 & 45.4806 \\
		\hline
		BASE          & 223.5281 & 223.5281 & 223.5281 & 223.5281 & 223.5281 \\
		\hline
	\end{tabular}
\end{table}

\begin{table}[t]
	\caption{Sensor network tracking problem: Mean Square Error (MSE) of various filtering schemes as a function of process noise parameter $\sigma_Q$. The boldfaces show the best performers for small/large particle sizes.}
	\label{table_sensor}
	\centering
	 \def\arraystretch{1.25}
	\begin{tabular}{ |c|c|c|c|c|c| } 
		\hline
		& \multicolumn{5}{|c|}{$\sigma_Q$} \\ \cline{2-6}
		&$10^{-2}$ & $5 \times 10^{-2}$ & $10^{-1}$ & $5 \times 10^{-1}$ & $1$ \\
		\hline
		SKF           & 10.4674 & 9.5812  & 9.5038  & 10.1664 & 10.5996 \\
		MKF & 10.5572 & 9.2879  & 9.1684  & 9.8175  & 10.3307 \\
		$\alpha$KF    & 9.9441 & 8.0913 & \textbf{8.0623} & 9.1002 & 9.7055 \\ 
		\hline
		PF            & 9.5661  & 9.3464  & 9.4726  & 10.0422 & 10.3834 \\ 
		\hline
		EKF           & 14.0034 & 13.9357 & 14.5161 & 15.5277 & 16.1438 \\
		UKF           & 11.5303 & 10.3639 & 10.2068 & 10.8830 & 11.5845 \\
		\hline
	\end{tabular}
\end{table}

\begin{figure*}[t]
	\centering
	\includegraphics[width=1\textwidth]{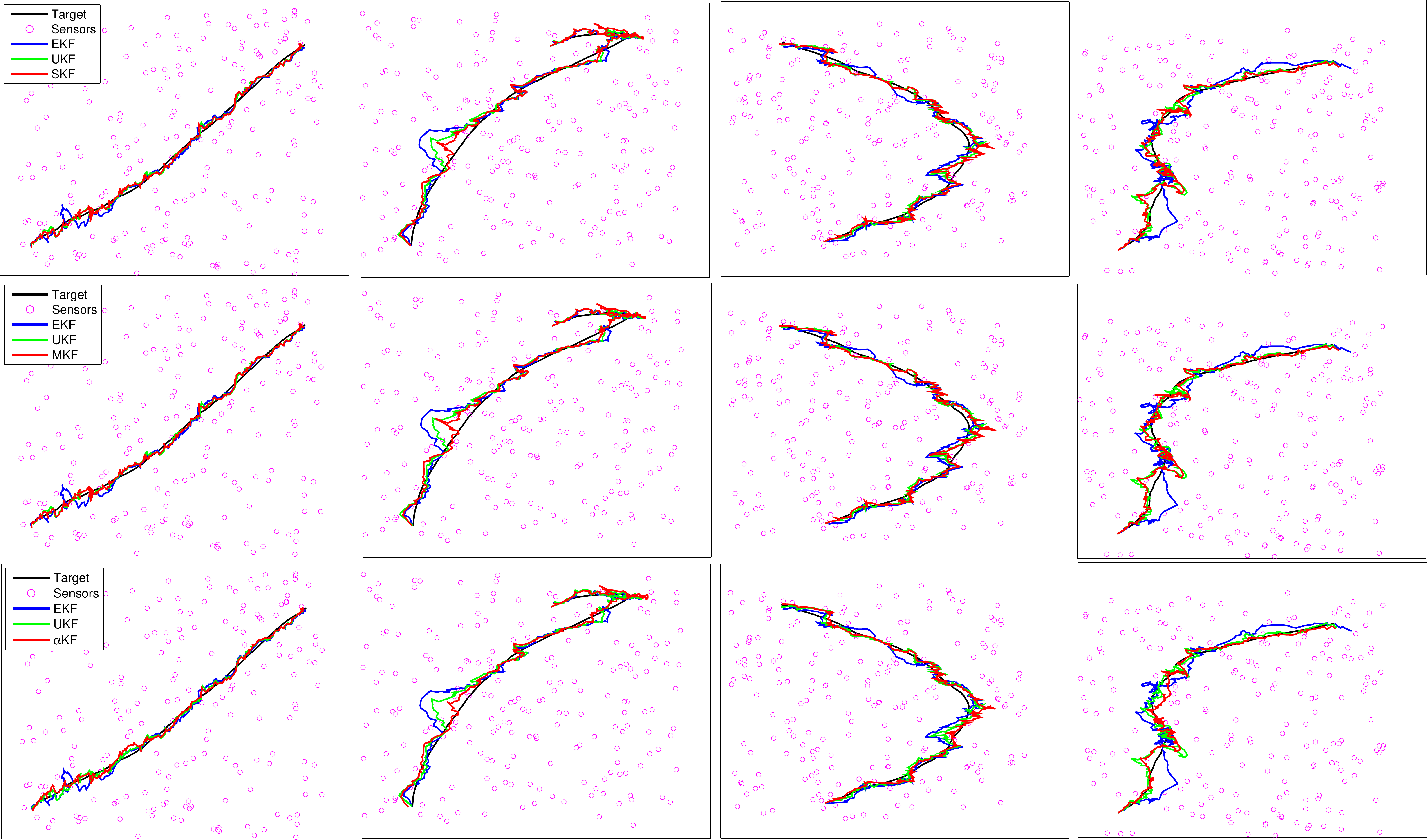}
	\caption{Tracks estimated by various filtering schemes in sensor network setting. Top row: Comparisons of EKF, UKF, and SKF. Middle row: EKF, UKF, and MKF. Bottom row: EKF, UKF, and {\akf} \hspace{-.07in}. Best viewed in color.}
	\label{fig:sens_traj}
\end{figure*}

In Table \ref{table_radar} we show mean square error (MSE) for radar tracking as a function of the selected scale value ($\sigma_Q$). Here, the base error corresponds to the estimations based on measurements only, and its order-of-magnitude difference from filter MSE values show the severity of nonlinearity. Now, comparing MSE values, first we see that MKF and \akf overperforms EKF and UKF for all settings of $\sigma_Q$ which shows that the Gaussian density obtained from these filters is indeed more accurate. SKF also gets better results, particularly for $\sigma_Q=10^{-1}$ but it is less robust to the changes in scale value. This is due to the iterative gradient scheme employed by SKF, which could give worse results depending on parameter changes or covariance initializations. Since MKF/\akf are based on importance sampling, they do not exhibit the same sensitivity. As for PF, this algorithm also produces competitive results when $\sigma_Q=10^{-2}$; however its performance significantly deteriorates (even 
more than that of SKF) as $\sigma_Q$ increases, which shows the nonparametric inference of particle filtering is more sensitive to parameter uncertainty. We also mark the best overall MSE with boldfaces, which is given by \akf for $\sigma_Q = 10^{-1}$. Furthermore, \akf has the highest robustness to parameter changes, therefore it is a better candidate to choose when parameters are not known and measurements are very noisy, since the $\alpha$ coefficient has the capability of mitigating excess measurement noise, as discussed in Section \ref{sec:alpha}.
 
Table \ref{table_sensor} presents MSE results for sensor networks. Unlike the radar problem, all particle-based filters are better than EKF/UKF for all values of $\sigma_Q$. This reduced sensitivity is due to the reduced nonlinearity in the problem. The performance of  SKF, MKF, and PF are similar to each other, MKF being the favorable choice for most of the cases. On the other hand \akf is the best performer in all cases, and as $\sigma_Q$ increases, the margin increases. The best overall MSE is again achieved by this filter for $\sigma_Q=10^{-1}$, where using \akf provides a clear benefit.
 
In Figure \ref{fig:sens_traj} we show qualitative tracking results from sensor networks. The top, middle, and bottom rows correspond to SKF, MKF, and \akf respectively. For each two we pick four different paths (shared across different rows) and for each plot we plot the true trajectory along with EKF, UKF, and one of our filters, depending on the row. First we see that our simulation settings encompass a wide variety of paths which exhibit multimodality such as, for example, a combination of constant velocity and constant turn models \cite{LiJilkov5}. By visual inspection we can see that our algorithms provide more accurate tracking compared to EKF/UKF in all cases. Furthermore, moving down the rows we can see that the accuracy of our filtering algorithms also increase and the \akf estimated paths are more robust to measurement errors, as clearly demonstrated in second column.

\begin{figure}[t]
\centering
	\includegraphics[width=.9\columnwidth]{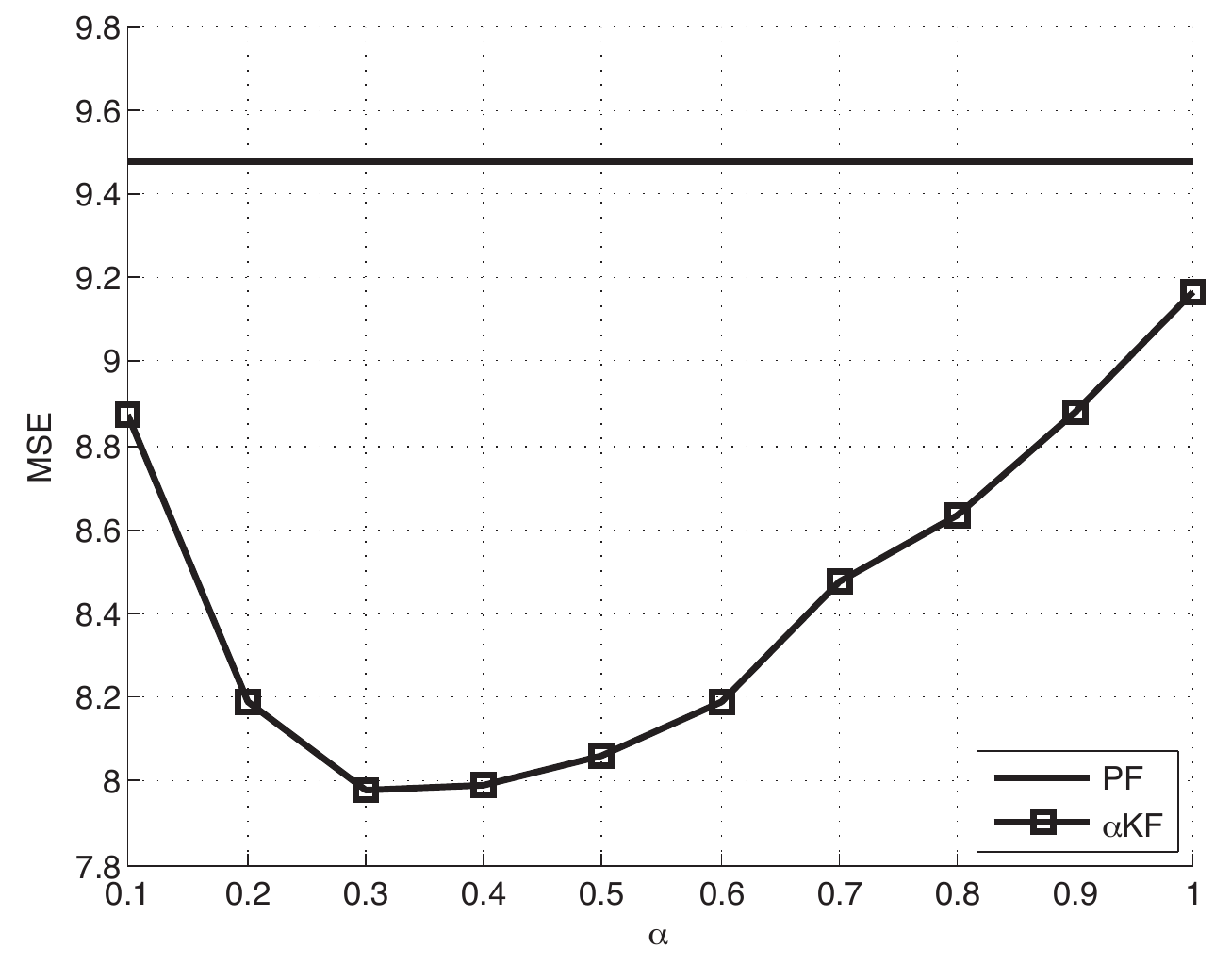}
	\caption{MSE value of \akf as a function of $\alpha$ for the sensor network tracking problem with $\sigma_Q=10^{-1}$. When $\alpha=1$, {\akf} reduces to MKF. The performance of PF is plotted as a baseline.}
	\label{fig:aplot}
\end{figure}

\begin{figure*}[t]
\centering
\subfigure[MSE as a function of $r_\text{max}$.]{\includegraphics[height=1.42in]{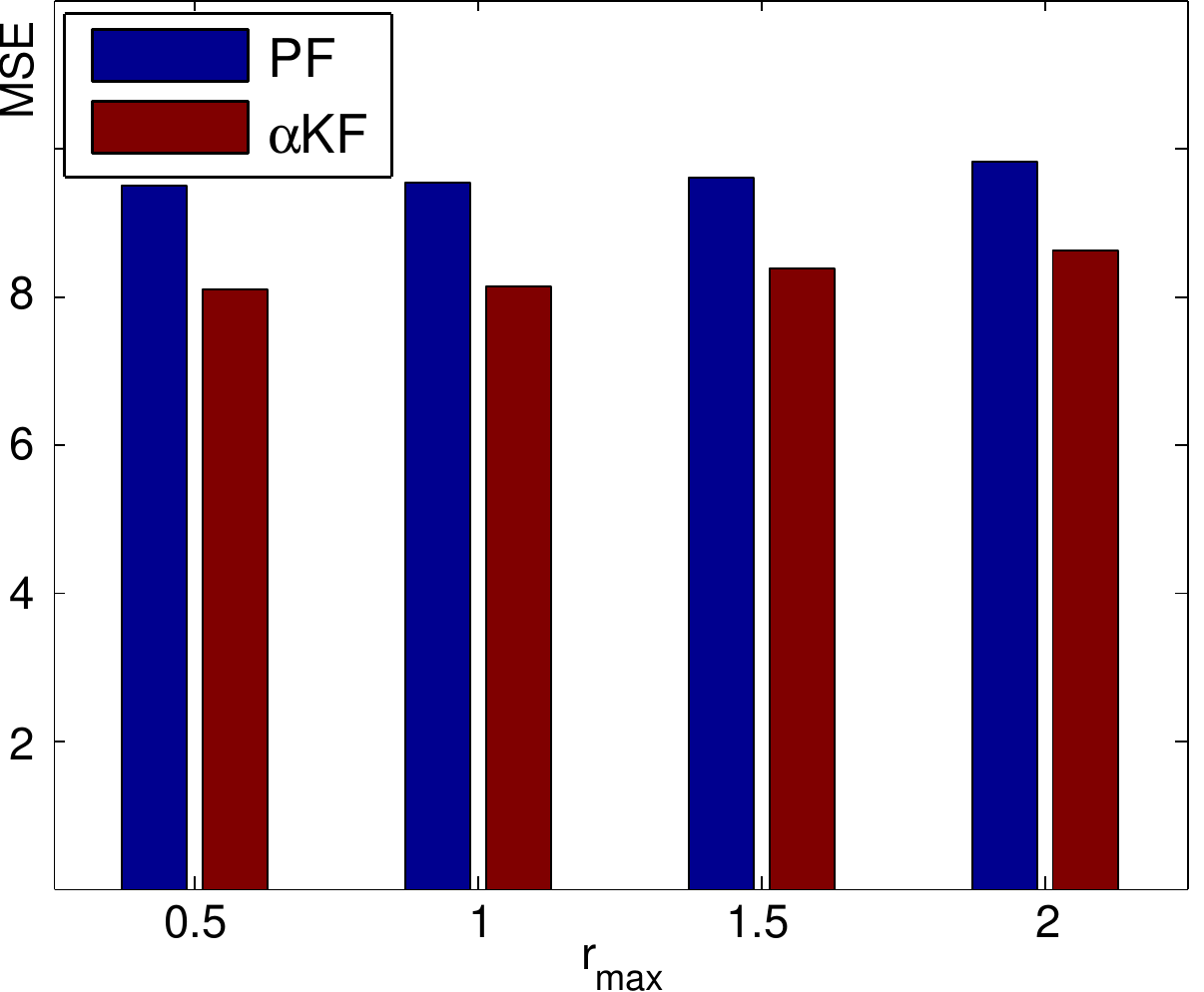}}
\subfigure[Samples required to achieve $r_\text{max}$.]{\includegraphics[height=1.42in]{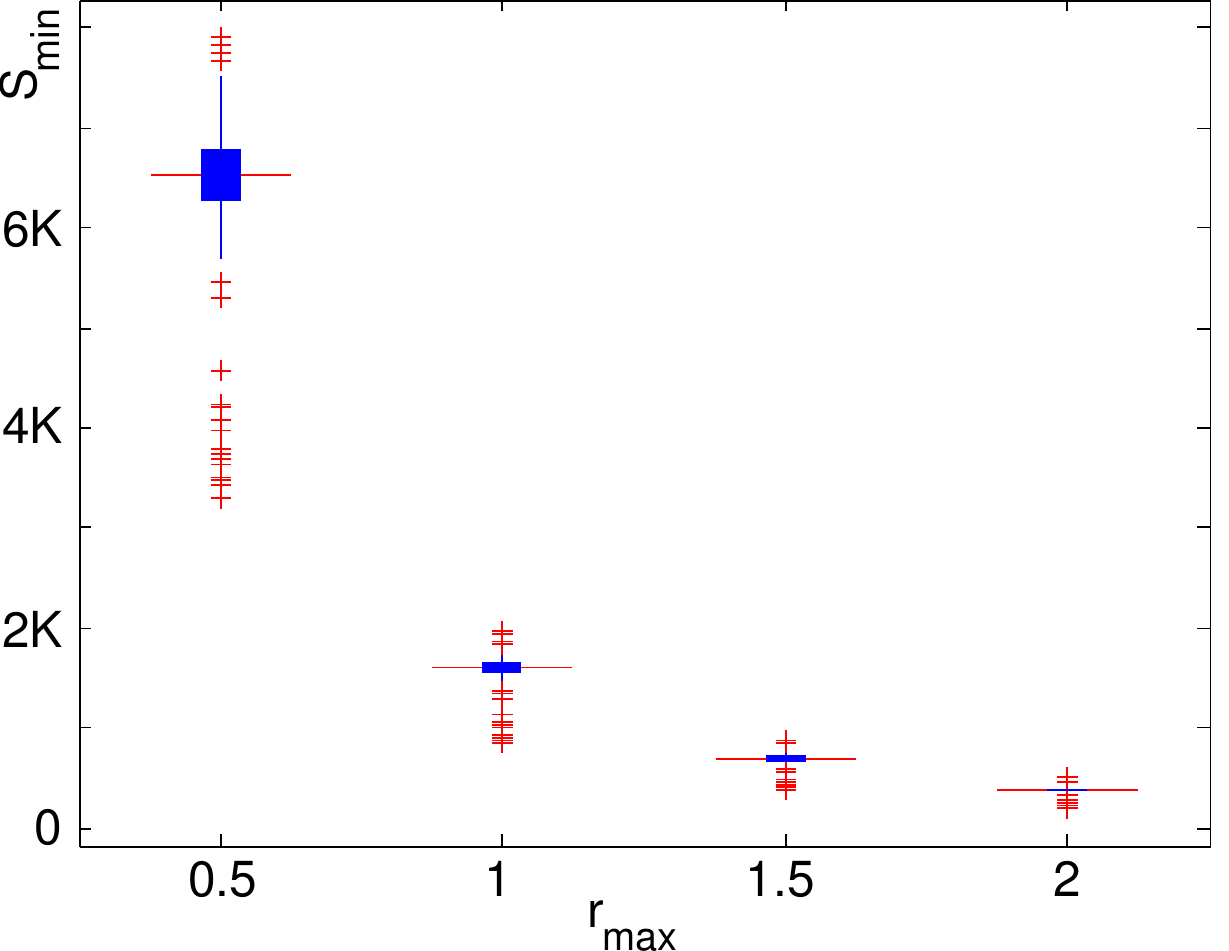}}
\subfigure[MSE as function of process noise.]{\includegraphics[height=1.42in]{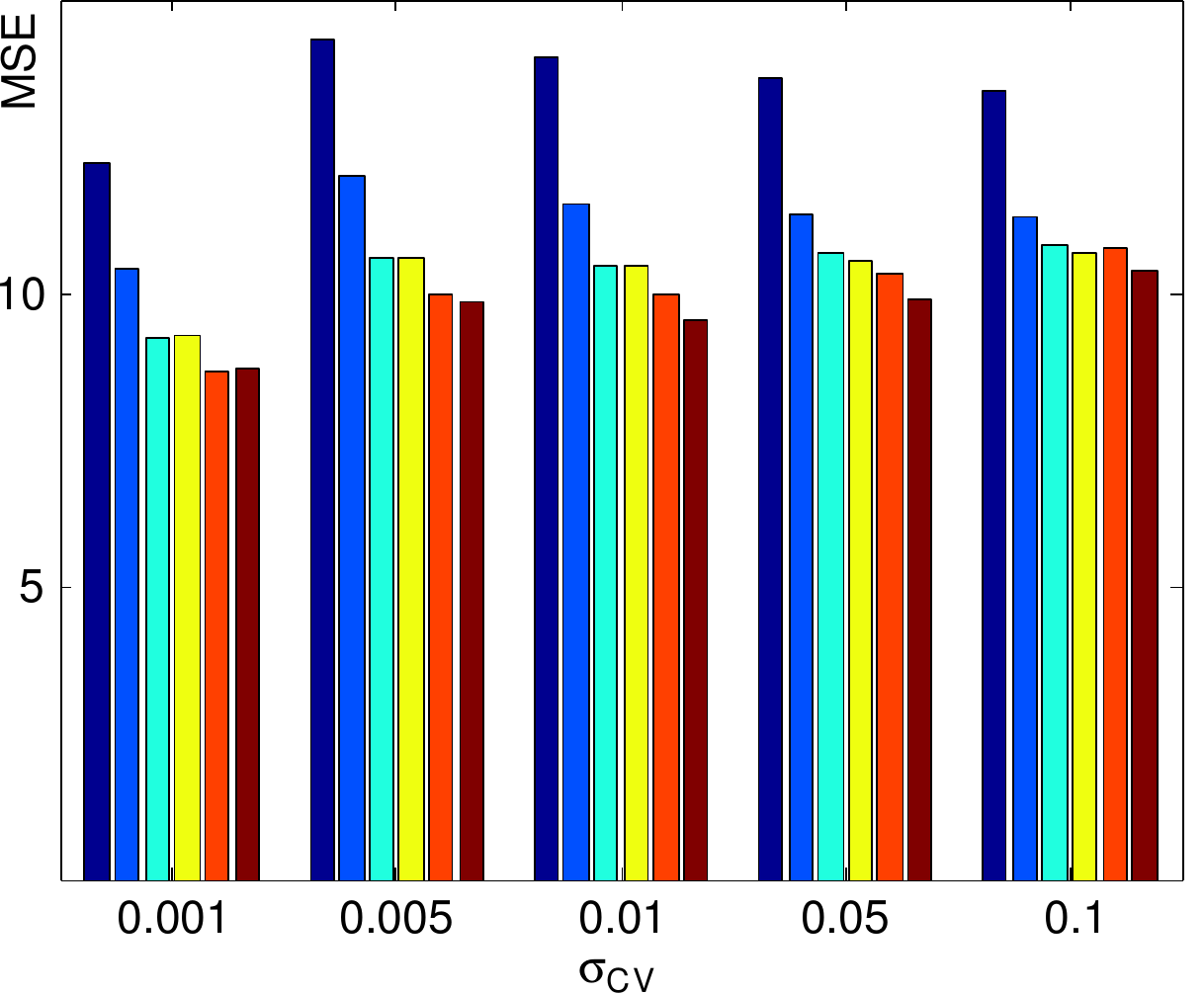}}
\subfigure[MSE versus measurement noise.]{\includegraphics[height=1.42in]{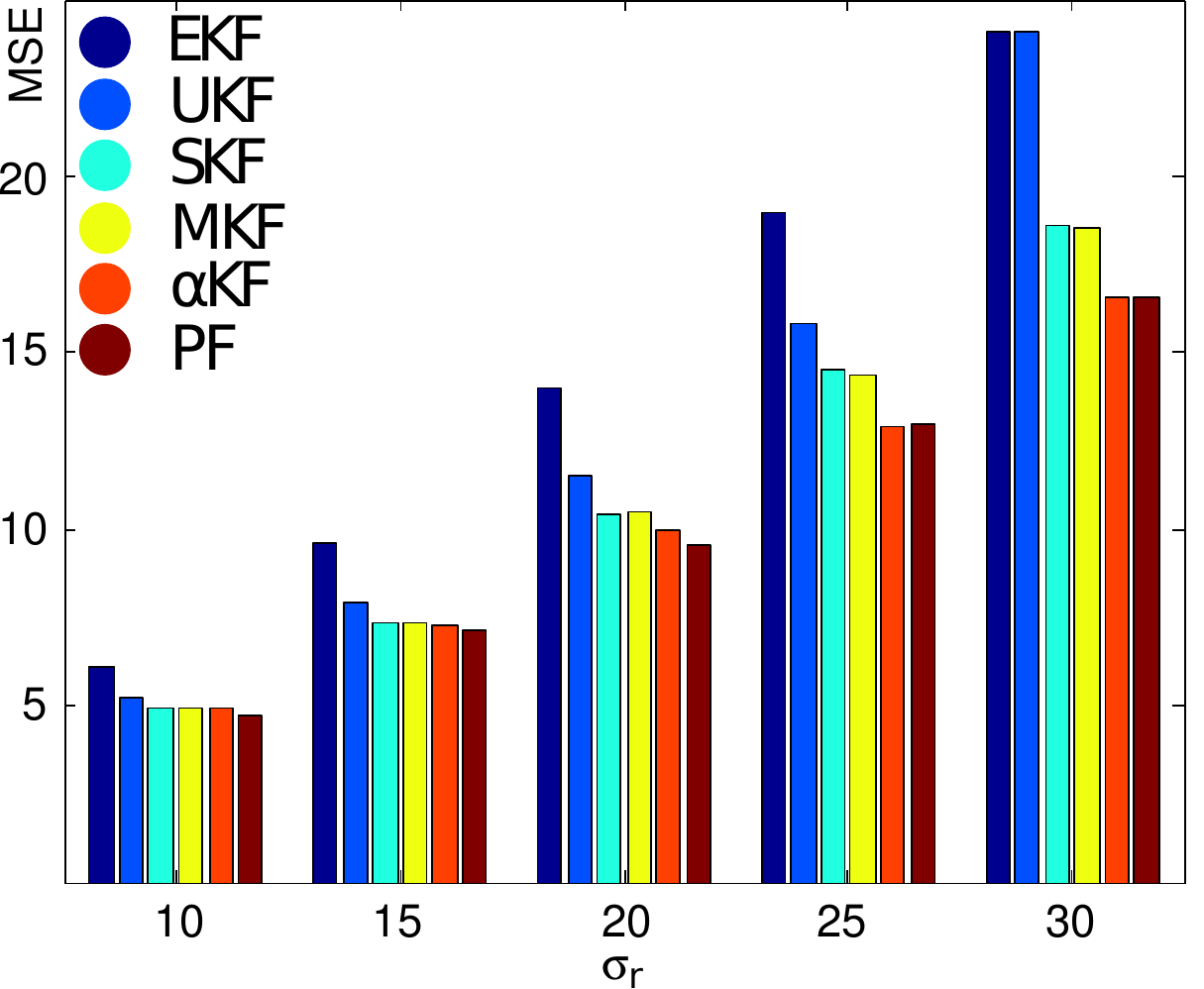}}
\caption{Mean square error and minimum sample size as a function of confidence radius $r_\text{max}$. Mean square error as a function of process and measurement noise parameters, where the exact parameters are known to the filter. The legend given is shared by both figures. (Best viewed in color.)}\label{fig:rmax_plot-qrplot}
\end{figure*}

So far, for \akf we only considered the case when $\alpha=0.5$, which used the symmetric Hellinger distance metric, as given in \eqref{eq:alp_cases}. Now we focus on varying the value of $\alpha$ and analyzing its effects. For this we use the sensor experiments with $\sigma_Q = 10^{-1}$, which corresponds to the mid column of Table \ref{table_sensor}. The mean squared error as a function of $\alpha$ is plotted in Figure \ref{fig:aplot}. We see that low-mid ranges of $\alpha$ (i.e $0.3-0.5$) give the best MSE results. This improvement is obtained since lower values of $\alpha$ help mitigate the effects of strong measurement noise. There is, however, a tradeoff here since choosing a too small value for this parameter will discard all the measurement information and give poor results. This is seen for lower values of $\alpha$, where decreasing the parameter degrades performance.

As discussed in Section \ref{sec:adasamp} we can use adaptive sampling to choose the minimum possible sample size to achieve a certain confidence region radius, $r_\text{max}$. We implemented adaptive sampling for {\akf} using an initial batch size of $S_\text{base} = 500$. We picked four different values of $r_\text{max}$ from $\{0.5,1,1.5,2\}$. Figure \ref{fig:rmax_plot-qrplot} displays the results for this experiment. In the left panel we compare the MSE results as a function of $r_\text{max}$ for {\akf} and PF for the sensor tracking problem with $\sigma_\text{CV}=0.1$. Note that, for PF, adaptive sampling is not a choice as all particles should be propagated, resampled, and updated at every time step. So for PF we simply set the sample size as the average $S_\text{min}$ for the {\akf} for each case. We can see that, the MSE performances differ very little across different cases, showing even for larger target values of $r_\text{max}$ both methods can still produce accurate estimates of the true state. We 
also 
see that {\akf} overperforms PF in all cases. On the other hand, the right panel shows the number of samples required to achieve a certain confidence radius. From this figure we can see the $\mathcal{O}(1/r^2)$ decaying rate of $S_\text{min}$ as implied by \eqref{eq:rmax_ratio}. Given the high accuracies in the left panel, we see that several hundred samples can be sufficient to obtain high-quality estimates, which makes {\akf} competitive for real time applications. Another point is, as $r_{max}$ increases, the variance of the sample size also decreases, which means the runtime per round will have small discrepancy, as opposed to using a smaller $r_\text{max}$.


We now turn to the case where the process noise parameter is known. In Figure \ref{fig:rmax_plot-qrplot} we show the filter MSEs as a function of $\sigma_{CV}$ and $\sigma_R$. For the measurement noise, as $\sigma_R$ increases the overall MSEs also increase, while for process noise this trend is not present. For both cases we see that the particle filter gives the best result overall. This is expected, since when the parameters are known perfectly, particle filter can approximate the posterior with more accuracy, as it is nonparametric. With that said, {\akf} is also competitive in this setting. In fact, for several cases such as $\sigma_\text{CV}=0.001$ and $\sigma_r=25$ performance of {\akf} and PF are equal, and for the remaining cases the particle filter does not improve much compared to {\akf}, while both filters can perform much better than SKF and MKF. This means {\akf} can be preferred over PF, since it does not require resampling. As a second observation, note that SKF/MKF perform much better than 
EKF/UKF, 
and {\akf} perform even better compared to the rest. This means, by minimizing different forms of divergence one can indeed get significantly better Gaussian approximations of the posterior, which supports our theoretical analysis in Section \ref{sec:new_nonkf}.



\subsection{Options Pricing}

We also consider a problem in options pricing. In finance, an option is a derivative security which gives the holder a right to buy/sell (call/put option) the underlying asset at a certain price on or before a specific date. The underlying asset can be, for example, a stock. The price and date are called the strike price and expiry date respectively. The value of the option, called premium, depends on a number of factors. Let $C$ and $P$ denote the call and put prices. We use $\sigma$ and $r$ to denote volatility and risk-free interest rate respectively; the values of these variables are not directly observed, hence they need to be estimated. Let $S$ denote the price of underlying asset and $X$ denote the strike price. Finally, let $t_m$ denote the time to maturity; this is the time difference between the purchase and expiry dates which is written as a fraction of a year. For example, an option which expires in two months will have $t_m = 1/6$.

Accurate pricing of options is an important problem in mathematical finance. For a European style option, the price as a function of all these parameters can be modeled using the well-known Black-Scholes equation \cite{Hull_2006}
\begin{align}
d_1 &= \frac{\log(S/X) + (r +\sigma^2/2)t_m}{\sigma\sqrt{t_m}} ~,~  d_2= d_1 - \sigma \sqrt{t_m} ~,~ \nonumber \\
C &= S \Phi(d_1) - X e^{-r t_m} \Phi(d_2) ~,~ \nn \\
P &= -S \Phi(-d_1) + X e^{-r t_m} \Phi(-d_2) \label{eq:call_put_price} ~.
\end{align}
Following the approach of \cite{Niranjan_1997}, let $x = {[\sigma ~ r]}^{\top}$ be the state and $y = {[C ~ P]}^{\top}$ be the measurement. We get the following state space representation
\begin{eqnarray}
x_t &= x_{t-1} + w_t ~,&~ w_t \sim N(0,Q)~, \nonumber \\
y_t &= h(x_t) + v_t ~,&~ v_t \sim N(0,R)~. \label{eq:options_model}
\end{eqnarray}
where the nonlinear mapping $h(\cdot)$ is given by \eqref{eq:call_put_price}. In this case we model the process and measurement noises with time-invariant covariance matrices $Q$ and $R$. We consider two tasks: 1) predicting the one-step ahead prices, and 2) estimating the values of hidden state variables. This problem is also considered in \cite{Merwe_2000} to assess the performance of particle filtering algorithms.

Here we use the Black-Scholes model as the ground truth. In order to synthesize the data, we use historical values of VIX (CBOEINDEX:VIX), which measures the volatility of S\&P 500 companies. From this list we pick Microsoft (NASDAQ:MSFT), Apple (NASDAQ:AAPL), and IBM (NYSE:IBM) as underlying assets and use their historical prices. The interest rate comes from a state-space model with a process noise of zero mean and variance $10^{-4}$. We set $\sigma_Q=\sigma_R=10^{-2}$. In Table \ref{table_option} we show the next-day prediction performance of all algorithms. We can see that the prediction performance imporves as we move towards MKF. This, again shows the difference between Gaussian approximations of the methods we employ. For MKF and PF we used $10^3$ particles, and their results were similar so we only report MKF here; however we also note that MKF can achieve this performance without using resampling, and it can leverage adaptive sampling to reduce sample size, which makes it preferable over PF. On the 
other hand, for SKF we need to use a large number of particles per iterations (around $1,000$). Even though this gives better results then EKF and UKF it is much slower than MKF/PF, and its performance can vary significantly between iterations, which makes it less competitive in this case. On the other hand, since the measurement noise is small in this case, choosing $\alpha < 1$ for \akf does not provide improvement over MKF in this case, which is consistent with our previous intuition. Therefore $\alpha=1$ is the best choice in this case.

\begin{figure}[t]
\centering
\includegraphics[width=1\columnwidth]{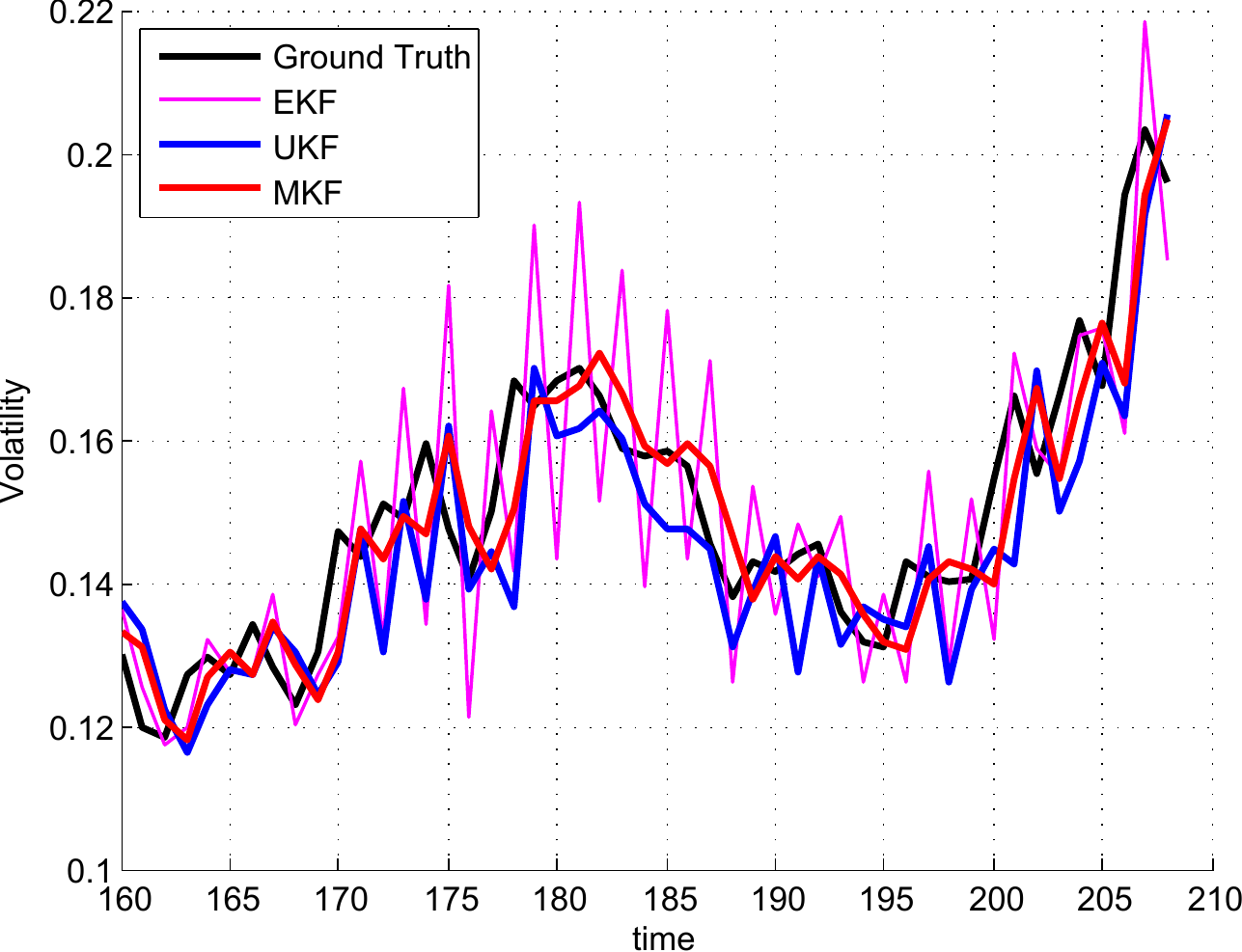}
\caption{Volatility estimation performance of various filtering schemes (based on Option 1). The estimates are plotted along with the ground truth. Best viewed in color.}
\label{fig:vol_perf}
\end{figure}

Figure \ref{fig:vol_perf} shows the volatility estimation for three filters: Usually EKF tends to over/under-shoot a lot and UKF is significantly better in that respect; however MKF improves even further as it gives the most robust estimates. The plot of SKF output is similar to MKF. Also, similar to the target tracking experiments, we see that MKF has better performance than SKF, which once again agrees with the observation that expectation-propagation typically outperforms variational inference for unimodal posterior.

\begin{table}[t]
	\caption{Mean Absolute Error(MAE) values of various filtering schemes for three different call/put option pairs; calculated for $\sigma_Q = 10^{-2}$. For Option 3, EKF loses track so MAE is not reported.}
	\label{table_option}
	\centering
	 \def\arraystretch{1.25}
	\begin{tabular}{ |c|c|c|c|c|c|c|c| } 
		\hline
		& & EKF & UKF & SKF & MKF \\
		\hline
		\multirow{2}{4em}{Option 1 \\ MAE} & Call & 0.1352 & 0.0788 & 0.0658 & 0.0654 \\
		& Put & 0.1528 & 0.0789 & 0.0642 & 0.0654 \\
		\hline
		\multirow{2}{4em}{Option 2 \\ MAE} & Call & 0.0425 & 0.0354 & 0.0312 & 0.0319 \\
		& Put & 0.0478 & 0.0355 & 0.0368 & 0.0331 \\
		\hline
		\multirow{2}{4em}{Option 3 \\ MAE} & Call & - & 0.2155 & 0.1573 & 0.1586 \\
		& Put & - & 0.2158 & 0.1574 & 0.1586 \\
		\hline
	\end{tabular}
\end{table}

\appendices

\section{ADF Equations} \label{app1}

For the proofs in the following appendices we need the predict-update equations of the joint Gaussian ADFs. Note that this corresponds to the model in \eqref{eq:jnt_asmp}. The equations are summarized as
\begin{eqnarray}
\text{Predict:} \quad &x_{t|t-1} &= F_t x_{t-1|t-1} ~, \nonumber \\
& P_{t|t-1} &= F_t P_{t-1|t-1} F_t^{\top} + Q_t ~, \label{eq:adf_prd} \\
\text{Update:} \quad &x_{t|t} &= x_{t|t-1} + K_t(y_t - y_{t|t-1}) ~, \nonumber \\
&P_{t|t} &= P_{t|t-1} - K_t S_t K_t^{\intercal} \label{eq:adf_upd} ~, \\
\text{Auxiliary:} \quad & y_{t|t-1} &= \mu_y = h_t(x_{t|t-1}) ~, \nonumber \\
& H_t &= \rSigma_{yx} \rSigma_{xx}^{-1} ~, \nonumber \\
& S_t &= \rSigma_{yy} = \rSigma_{yx}\rSigma_{xx}^{-1}\rSigma_{xy} + R_t ~, \nonumber \\
& K_t &= \rSigma_{xy} \rSigma_{yy}^{-1} ~. \label{eq:adf_aux}
\end{eqnarray} 
We emphasize that these hold for \emph{any} joint Gaussian ADF. When EKF is employed, $H_t$ is the Jacobian at prior mean, and $S_t$ and $K_t$ are calculated accordingly.

\section{Proof of Theorem 1} \label{app2}

The joint Gaussian ADF corresponds to $f(x) \approx g(x)$; this approximation is constructed from $p(y_t|x_t)$ in \eqref{eq:jnt_asmp}, which is Gaussian with $\mu_{y|x} = \mu_y + \rSigma_{yx} \rSigma_{xx}^{-1} (x_t-\mu_{x})$ and $\rSigma_{y|x} = \rSigma_{yy} - \rSigma_{yx} \rSigma_{xx}^{-1} \rSigma_{xy}$. This yields
\begin{equation} 
	g(x_t) = -\frac{1}{2}(\wty_t - \rSigma_{yx} \rSigma_{xx}^{-1} x_t )^{\top} R_t^{-1} (\wty_t - \rSigma_{yx} \rSigma_{xx}^{-1} x_t ) ~, \label{eq:adf_g}
\end{equation}
where $\wty_t = y_t - \mu_y + \rSigma_{yx} \rSigma_{xx}^{-1} \mu_{x}$. Note that under \eqref{eq:jnt_asmp} we have $p(x_t) \sim N(\mu_x,\rSigma_{xx})$, and let $q(x_t) \sim N(\hmu,\hSigma)$. Substituting $g(x_t)$ to \eqref{eq:elbo_kf} the expectations are now evaluated as
\begin{align}
	-\mathbb{E}_q [\log q(x_t)] &= \frac{1}{2} \log |\hSigma_t| ~, \nn \\
	-\mathbb{E}_q [\log p(x_t)] &= -\frac{1}{2} \hmu^{\top} \rSigma_{xx}^{-1} \hmu - \frac{1}{2} \tr\{\rSigma_{xx}^{-1} \hSigma_t\} + \hmu_t \rSigma^{-1} \mu_x ~, \nn \\
	-\frac{1}{2} \mathbb{E}_q[g(x_t)] &= -\frac{1}{2} \hmu_t^{\top} \rSigma_{xx}^{-1} \rSigma_{xy} R_t^{-1} \rSigma_{yx} \rSigma_{xx}^{-1} \hmu_t \nn \\ 
	&\quad -\frac{1}{2} \tr \{ (\rSigma_{xx}^{-1} \rSigma_{xy} R_t^{-1} \rSigma_{yx} \rSigma_{xx}^{-1}) \hSigma_t \} \\ 
	&\quad +\hmu_t \rSigma_{xx}^{-1} \rSigma_{xy} R_t^{-1} \wty_t ~. \label{eq:adf_exp}	
\end{align}
The posterior parameters are found by solving $\nabla_{\hmu}\mathcal{L}=0$ and $\nabla_{\hSigma}\mathcal{L}=0$. Differentiating the terms in \eqref{eq:adf_exp} we get
\begin{align}
	\hSigma_t &= [\rSigma_{xx}^{-1} + \rSigma_{xx}^{-1} \rSigma_{xy} R_t^{-1} \rSigma_{yx} \rSigma_{xx}^{-1}]^{-1} ~, \label{eq:adf_sig1} \\
	\hmu_t &= \hSigma_t (\rSigma_{xx}^{-1}\mu_x + \rSigma_{xx}^{-1} \rSigma_{xy} R^{-1} \wty_t)	 ~. \label{eq:adf_mu1}
\end{align}
The matrix inversion lemma asserts $(A+UCV)^{-1} = A^{-1} - A^{-1}U(C^{-1}+VA^{-1} U)^{-1}VA^{-1}$; applying this to \eqref{eq:adf_sig1} we obtain
\begin{align}
	\hSigma_t &= \rSigma_{xx} - \rSigma_{xy} (\rSigma_{yx}\rSigma_{xx}^{-1}\rSigma_{xy} + R_t^{-1}) \rSigma_{yx} ~, \nn \\
	& = \rSigma_{xx} - \rSigma_{xy} \rSigma_{yy}^{-1} \rSigma_{yy} \rSigma_{yy}^{-1} \rSigma_{yx} ~. \label{eq:adf_sig2}	
\end{align}
Substituting \eqref{eq:adf_sig2} into \eqref{eq:adf_mu1} and expanding we get
\begin{align}
	\hmu_t &= \mu_x - \rSigma_{xy} \rSigma_{yy}^{-1} \rSigma_{yx} \rSigma_{xx}^{-1} \mu_x \nn \\
	&\quad + \rSigma_{xy} R_t^{-1} \wty_t - \rSigma_{xy} \rSigma_{yy}^{-1} \rSigma_{yx} \rSigma_{xx}^{-1} \rSigma_{xy} R_t^{-1} \wty_t ~, \nn \\
	&= \mu_x - \rSigma_{xy} \rSigma_{yy}^{-1} \rSigma_{yx} \rSigma_{xx}^{-1} \mu_x \nn \\
	&\quad + (I-\rSigma_{xx} \rSigma_{yy}^{-1} \rSigma_{yx} \rSigma_{xx}^{-1}) \rSigma_{xy} R_t^{-1} \wty_t ~, \nn \\
	&= \mu_x - \rSigma_{xy} \rSigma_{yy}^{-1} \rSigma_{yx} \rSigma_{xx}^{-1} \mu_x + \rSigma_{xy} \rSigma_{yy}^{-1} \wty_t ~, \nn \\
	&= \mu_x + \rSigma_{xy}\rSigma_{yy}^{-1}(y_t-\mu_y) ~. \label{eq:adf_mu2}
\end{align}
Note the third line follows from the identity $\rSigma_{xy} \rSigma_{yy}^{-1} = (I - \rSigma_{xy} \rSigma_{yy}^{-1} \rSigma_{yx} \rSigma_{xx}^{-1}) \rSigma_{xy} R_t^{-1}$ which can be verified with straightforward manipulation. Matching the terms in \eqref{eq:adf_aux} with \eqref{eq:adf_mu2} and \eqref{eq:adf_sig2} we obtain the updates in \eqref{eq:adf_upd}.$\quad \blacksquare$

\section{Proof of Corollary 2} \label{app3}

The proof is similar to that of Theorem 1, therefore we highlight the key points. We simplify the notation to $p(x_t) \sim N(\mu_t,\rSigma_t)$ and $q(x_t) \sim N(\hmu_t,\hSigma_t)$. We employ a first-order Taylor series expansion around prior mean: $h(x_t) \approx h(\mu_t) + H_t(\mu_t)(x_t-\mu_t)$ where $H_t$ is the Jacobian. Define $\wty_t = y_t - h(\mu_t) + H_t\mu_t$. Plugging these into the variational lower bound \eqref{eq:elbo_kf} and differentiating we obtain 
\begin{align}
\hSigma_t &= {(\rSigma_t^{-1} + H_t^{\top}R_t^{-1}H_t)}^{-1} ~, \label{eq:ekf_sig1} \\
\hmu_t &= \hSigma_t ~ (\rSigma_t^{-1}\mu_t + H_t^{\top} R_t^{-1} \wty_t) ~. \label{eq:ekf_mu1} 
\end{align}
Once again, using the matrix inversion lemma we get
\begin{align}
\hSigma_t &= \rSigma_t - K_t S_t K_t^{\top} ~, \label{eq:ekf_sig2}
\end{align}
where $S_t = H_t \rSigma_t H_t^{\top} + R_t$ and $K_t = \rSigma_t H_t^{\top} {S_t}^{-1}$. Plugging \eqref{eq:ekf_sig2} in \eqref{eq:ekf_mu1} and expanding the multiplication we get 
\begin{align}
\hmu_t &= \mu_t + \rSigma_t H_t^{\top} R_t^{-1} \wty_t - K_t H_t^{\top} \mu_t - K_t H_t \rSigma_t H_t^{\top} R_t^{-1} \wty_t \nonumber \\
&= \mu_t - K_t H_t^{\top} \mu_t + (I -K_t H_t) \rSigma_t H_t^{\top} R_t^{-1} \gamma \nonumber \\
&= \mu_t  - K_t H_t^{\top} \mu + K_t \wty_t \nonumber \\
&= \mu_t + K_t (y_t - h_t(\mu_t))  \label{eq:ekf_mu2}
\end{align}
where the first and third lines utilize the identities $\hSigma_t = \rSigma_t - K_t H_t \rSigma_t$ and $K_t = (I-K_t H_t) \rSigma_t H_t^{\top}R_t^{-1}$ respectively. We see that \eqref{eq:ekf_mu2} and \eqref{eq:ekf_sig2} correspond to the EKF update equations. $\quad \blacksquare$

\bibliographystyle{IEEEbib}
\bibliography{nonkf_bib}

\begin{thebibliography}{10}

\bibitem{Kalman_1960}
Rudolph~Emil Kalman,
\newblock ``A new approach to linear filtering and prediction problems,''
\newblock {\em Transactions of the ASME--Journal of Basic Engineering}, vol.
  82, no. Series D, pp. 35--45, 1960.

\bibitem{Belanger_2015}
David Belanger and Sham Kakade,
\newblock ``A linear dynamical system model for text,''
\newblock in {\em International Conference on Machine Learning (ICML)}, 2015.

\bibitem{Gultekin_2014}
San Gultekin and John Paisley,
\newblock ``A collaborative kalman filter for time-evolving dyadic processes,''
\newblock in {\em 2014 IEEE International Conference on Data Mining}, Dec 2014,
  pp. 140--149.

\bibitem{Blei_2006}
David~M. Blei and John~D. Lafferty,
\newblock ``Dynamic topic models,''
\newblock in {\em ICML}, 2006.

\bibitem{LiJilkov3}
Xiao~Rong Li and Vesselin~P. Jilkov,
\newblock ``Survey of maneuvering target tracking: Iii. measurement models,''
\newblock in {\em International Symposium on Optical Science and Technology},
  2001.

\bibitem{Koren_2009}
Yehuda Koren, Robert Bell, and Chris Volinsky,
\newblock ``Matrix factorization techniques for recommender systems,''
\newblock {\em Computer}, vol. 42, no. 8, pp. 30--37, Aug 2009.

\bibitem{Welch_1995}
Greg Welch and Gary Bishop,
\newblock ``An introduction to the {K}alman filter,''
\newblock Tech. {R}ep., Chapel Hill, NC, US, 1995.

\bibitem{Julier_2004}
Simon~K. Julier and Jeffrey~K. Uhlmann,
\newblock ``Unscented filtering and nonlinear estimation,''
\newblock {\em Proceedings of the IEEE}, 2004.

\bibitem{Arulampalam_2002}
M.~Sanjeev Arulampalam, Simon Maskell, Neil Gordon, and Tim Clapp,
\newblock ``A tutorial on particle filters for online nonlinear/non-gaussian
  bayesian tracking,''
\newblock {\em IEEE Transactions on Signal Processing}, vol. 50, no. 2, pp.
  174--188, Feb 2002.

\bibitem{Jordan_1999}
Michael~I. Jordan, Zoubin Ghahramani, Tommi~S. Jaakkola, and Lawrence~K. Saul,
\newblock ``An introduction to variational methods for graphical models,''
\newblock {\em Machine Learning}, 1999.

\bibitem{Minka_2001}
Thomas~P. Minka,
\newblock ``Expectation propagation for approximate bayesian inference,''
\newblock in {\em Uncertainty in Artificial Intelligence (UAI)}, 2001.

\bibitem{Beal_2003}
Matthew~J. Beal,
\newblock {\em Variational Algorithms for Approximate Bayesian Inference},
\newblock Ph.D. thesis, University of London, 2003.

\bibitem{Vermaak_2003}
J.~Vermaak, N.~Lawrence, and P.~Perez,
\newblock ``Variational inference for visual tracking,''
\newblock in {\em Computer Vision and Pattern Recognition (CVPR)}, 2003.

\bibitem{Snoussi_2012}
Jing Teng, Hichem Snoussi, Cedric Richard, and Rong Zhou,
\newblock ``Distributed variational filtering for simultaneous sensor
  localization and target tracking in wireless sensor networks,''
\newblock {\em IEEE Transactions on Vehicular Technology}, 2012.

\bibitem{Maybeck_1982}
P.~S. Maybeck,
\newblock {\em Stochastic models, estimation, and control},
\newblock Academic Press Inc., 1982.

\bibitem{Ito_2000}
K.~Ito and K.~Xiong,
\newblock ``Gaussian filters for nonlinear filtering problems,''
\newblock {\em IEEE Transactions on Automatic Control}, 2000.

\bibitem{Guo_2006}
Dong Guo and Xiaodong Wang,
\newblock ``Quasi-monte carlo filtering in nonlinear dynamic systems,''
\newblock {\em IEEE Transactions on Signal Processing}, 2006.

\bibitem{Binjia_2013}
Bin Jia, Ming Xin, and Yang Cheng,
\newblock ``High-degree cubature kalman filter,''
\newblock {\em Automatica}, 2013.

\bibitem{Heskes_2002}
Tom Heskes and Onno Zoeter,
\newblock ``Expectation propagation for approximate inference in dynamic
  bayesian networks,''
\newblock in {\em Uncertainty in Artificial Infelligence (UAI)}, 2002.

\bibitem{Andrieu_2003}
Christophe Andrieu, Nando De~Freitas, Arnaud Doucet, and Michael~I. Jordan,
\newblock ``An introduction to mcmc for machine learning,''
\newblock {\em Machine learning}, vol. 50, no. 1-2, pp. 5--43, 2003.

\bibitem{Kotecha_2003}
Jayesh~H. Kotecha and Peter~M. Djuric,
\newblock ``Gaussian particle filtering,''
\newblock {\em IEEE Transactions on Signal Processing}, vol. 51, no. 10, pp.
  2592--2601, 2003.

\bibitem{Wainwright_2008}
Martin Wainwright and Michael Jordan,
\newblock ``Graphical models, exponential families, and variational
  inference,''
\newblock {\em Foundations and Trends in Machine Learning}, vol. 1, no. 1-2,
  2008.

\bibitem{Bishop_2006}
Christopher~M. Bishop,
\newblock {\em Pattern Recognition and Machine Learning},
\newblock Springer-Verlag New York, Inc., 2006.

\bibitem{Paisley_2012}
John Paisley, David~M. Blei, and Michael~I. Jordan,
\newblock ``Variational {B}ayesian inference with stochastic search,''
\newblock in {\em International Conference on Machine Learning (ICML)}, 2012.

\bibitem{Amari_1998}
Shun-Ichi Amari,
\newblock ``Natural gradient works efficiently in learning,''
\newblock {\em Neural Computation}, vol. 10, no. 2, pp. 251--276, Feb. 1998.

\bibitem{Yi_2009}
Sun Yi, Daan Wiestra, Tom Schaul, and Jurgen Schmidhuber,
\newblock ``Stochastic search using the natural gradient,''
\newblock in {\em International Conference on Machine Learning (ICML)}, 2009.

\bibitem{Lobato_2016}
J.~F.~M. Hernandez-Lobato, Yingzhen Li, Mark Rowland, Daniel Hernandez-Lobato,
  Thang~D. Bui, and Richard~E. Turner,
\newblock ``Black box alpha divergence minimization,''
\newblock in {\em International Conference on Machine Learning (ICML)}, 2016.

\bibitem{Darling_2015}
Jacob~E. Darling and Kyle~J. DeMars,
\newblock ``Minimization of the kullback-leibler divergence for nonlinear
  estimation,''
\newblock in {\em Proceedings of the Astrodynamics Specialist Conference},
  2015.

\bibitem{Minka_2004}
Thomas~P. Minka,
\newblock ``Power ep,''
\newblock Tech. {R}ep., 2004.

\bibitem{Chong_2003}
Chee~Y. Chong and Srikanta~P. Kumar,
\newblock ``Sensor networks: Evolution, opportunities, and challenges,''
\newblock {\em Proceedings of the IEEE}, 2003.

\bibitem{Tubaishat_2003}
Malik Tubaishat and Sanjay Madria,
\newblock ``Sensor networks: An overview,''
\newblock {\em IEEE Potentials}, 2003.

\bibitem{Boukerche_2007}
A.~Boukerche, H.A.B. Oliveira, E.F. Nakamura, and A.A.F. Loureiro,
\newblock ``Localization systems for wireless sensor networks,''
\newblock {\em IEEE Wireless Communications}, 2007.

\bibitem{LiJilkov5}
Xiao~Rong Li and Vesselin~P. Jilkov,
\newblock ``Survey of maneuvering target tracking. part v. multiple model
  methods,''
\newblock {\em IEEE Transactions on Aerospace and Electronic Systems}, 2005.

\bibitem{LiJilkov1}
Xiao~Rong Li and Vesselin~P. Jilkov,
\newblock ``Survey of maneuvering target tracking. part i. dynamic models,''
\newblock {\em IEEE Transactions on Aerospace and Electronic Systems}, 2003.

\bibitem{Hull_2006}
John~C. Hull,
\newblock {\em Options, futures, and other derivatives},
\newblock Pearson, Prentice Hall, 2006.

\bibitem{Niranjan_1997}
Mahesan Niranjan,
\newblock ``Sequential tracking in pricing financial options using model based
  and neural network approaches,''
\newblock in {\em Neural Information Processing Systems (NIPS)}, 1997.

\bibitem{Merwe_2000}
Rudolph van~der Merwe, Arnaud Doucet, Nando de~Freitas, and Eric Wan,
\newblock ``The unscented particle filter,''
\newblock in {\em Neural Information Processing Systems (NIPS)}, 2000.

\end{thebibliography}

\end{document}